\documentclass[11pt,letterpaper]{amsart}

\usepackage{color}

\DeclareMathAlphabet{\mathcal}{OMS}{cmsy}{m}{n}

\usepackage{graphicx, epstopdf}
\usepackage{cases}
\usepackage{bbm}
\usepackage{amssymb}
\usepackage{txfonts}
\usepackage{amscd}
\usepackage{amsfonts,latexsym,amsmath,amsxtra,mathdots,amssymb,latexsym,mathtools}
\usepackage{mathabx}
\usepackage[all,cmtip]{xy}

\usepackage{colordvi}
\usepackage{multicol}
\usepackage{hyperref}
\usepackage{tikz}
\usepackage{float}
\usepackage{setspace, floatflt}

\usepackage{tensor}

\usepackage[normalem]{ulem}

\allowdisplaybreaks

\newcommand{\BA}{{\mathbb {A}}}  \newcommand{\CF}{{\mathcal {F}}} 
\newcommand{\BC}{{\mathbb {C}}} \newcommand{\CE}{{\mathcal {E}}} 
 \newcommand{\CA}{{\mathcal {A}}} 
\newcommand{\BG}{{\mathbb {G}}} \newcommand{\CU}{\mathcal{U}}

 \newcommand{\BR}{{\mathbb {R}}}

 \newcommand{\BZ}{{\mathbb {Z}}}

\newcommand{\GL}{{\mathrm {GL}}}

 \newcommand{\Tr}{{\mathrm{Tr}}}

\newcommand{\res}{\mathrm{Res}}

\newcommand{\Hom}{\mathrm{Hom}}
\newcommand{\ind}{\mathrm{ind}} \newcommand{\Ind}{\mathrm{Ind}} \DeclareMathOperator{\cind}{c-ind}

\newcommand{\ra}{\rightarrow}

\def\-{^{-1}}

\renewcommand{\Re}{{\mathrm{Re}\,}}

\def\shskip{\hskip 0.5 pt}

\makeatletter
\g@addto@macro\normalsize{\setlength\abovedisplayskip{3pt}}
\makeatother

\makeatletter
\g@addto@macro\normalsize{\setlength\belowdisplayskip{3pt}}
\makeatother

\newcommand{\delete}[1]{}

\theoremstyle{plain}

\newtheorem{thm}{Theorem}[section] 
\newtheorem{lem}[thm]{Lemma}  \newtheorem{prop}[thm]{Proposition}

 \newtheorem{defn}[thm]{Definition}

\newtheorem {rem}[thm]{Remark}

\numberwithin{equation}{section}

\newcommand{\period}{\mathbf{P}}
\newcommand{\speh}{\mathrm{Sp}}
\newcommand{\wo}{w_{\star}}\newcommand{\CP}{\mathcal{P}}\newcommand{\CS}{\mathcal{S}}
\newcommand{\rat}{\textup{Rat}}

\begin{document}

	\title{Residual spectrum of $\GL_{2n}$ distinguished by $\GL_n \times \GL_n$}

	\author{Chang Yang}
	\address{Key Laboratory of High Performance Computing and Stochastic Information Processing (HPCSIP)\\ College of Mathematics and Statistics \\ Hunan Normal University \\Changsha,  410081\\China}
	\email{cyang@hunnu.edu.cn}


	\begin{abstract}
		Following the regularization method presented by Zydor, we study in this paper the regularized linear periods of square-integrable automormphic forms on $\GL_{2n}(\BA_F)$, where $F$ is a number field and $\BA_F$ its ring of adeles. We obtain a formula that expresses the regularized period of a noncuspidal, square-integrable automorphic form in terms of degenerate Whittaker functions in an inductive manner. As a consequence we characterize irreducible automorphic representations in the discrete spectrum of $\GL_{2n}(\BA)$ that are distinguished by $\GL_n(\BA) \times \GL_n(\BA)$. We also show the vanishing of the regularized periods of square-integrable automorphic forms on $\GL_n(\BA)$ over $\GL_p(\BA) \times \GL_q(\BA)$ when $p$ is not equal to $q$.
	\end{abstract}
	
	\maketitle

	\section{Introduction}
	
	Let $F$ be a number field and $\BA$ be its ring of adeles. Let $G$ be a reductive algebraic group over $F$ and $G'$ be a closed subgroup of $G$ defined over $F$. The period integral is defined as
	\begin{align}\label{formula::intro--denf-period}
		 \int_{G'(F) \backslash G'(\BA)} \varphi(g) \,dg
	\end{align}
	for an automorphic form $\varphi$ on $G(\BA)$, whenever it converges. Automorphic representations of $G(\BA)$ on which the period integral \eqref{formula::intro--denf-period} is not identically zero are called $G'$-distinguished. In many cases, distinguished representations can be characterized by functoriality and period integrals are closely related to special values of $L$-functions. Distinguished representations are those that contribute to the spectral decomposition of the theta series on $G/G'(\BA)$. To get a better understanding of these spectral decomposition of theta series, it is reasonable to determine distinguished representations in the entire automorphic spectrum of $G(\BA)$, not only the cuspidal spectrum.
	
	When $\varphi$ is not a cuspidal form, the integral \eqref{formula::intro--denf-period} may diverge. To regularize the period integral, Jacquet, Lapid and Rogawski introduced in \cite{J-L-R-periods-JAMS} a mixed truncation operator $\Lambda_m^T$ for split Galois pairs that is a relative version of Arthur's truncation operator. The mixed truncation operator was later generalized to all Galois pairs by Lapid and Rogawski \cite{Lapid-Rogawski-periods-Galois}, to symplectic pairs by Offen \cite{Offen-symplectic-disc-spectrum-IsraelJournal} and to the pair $(\GL_{n+1},\GL_n)$ by \cite{Ichino-Yamana-Compositio-Periods}. In a recent beautiful paper \cite{Zydor-periods-AENS}, the construction was generalized greatly by Zydor to all pairs $(G,G')$ where $G'$ is a reductive subgroup of $G$.
	
	In many cases, the period integral \eqref{formula::intro--denf-period} is closely related to a partial $L$-function via a Rankin-Selberg integral presentation of the $L$-function. The theory typically involves an integral of the form
	\begin{align}\label{formula::intro--Rankin-selberg-integral}
		\int_{G'(F) \backslash G'(\BA)} \varphi(g) E(g,\Phi,s ) dg,
	\end{align}
	where $\varphi$ is a cuspidal automorphic form on $G(\BA)$ and $E(\Phi,s)$ is an Eisenstein series on $G'(\BA)$. On one hand, the Eisenstein series $E(\Phi,s)$ often has a constant residue at some point $s = s_0$. Taking the residue at $s = s_0$, we then get the period integral. On the other hand, the integral \eqref{formula::intro--Rankin-selberg-integral} can be unfolded to a product of local zeta integrals that give local $L$-functions at almost all places of $F$. 
	
	In this paper we study regularized period integrals of noncuspidal automorphic forms by extending the above Rankin-Selberg method. This strategy was first carried out for the Galois pair $(\GL_n(E),\GL_n(F))$ with $E$ a quadratic field extension of $F$ by Yamana in \cite{Yamana-Residue-Periods-JFA} . The integral \eqref{formula::intro--Rankin-selberg-integral} may also diverge when $\varphi$ is not cuspidal. To regularize \eqref{formula::intro--Rankin-selberg-integral} amounts to regularize the period integral for the pair $(G \times G',\Delta G')$, where $\Delta G'$ is the diagonal embedding of $G'$ into $G \times G'$. Regularization procedure of such kind have been considered in \cite{Ichino-Yamana-Compositio-Periods} and \cite{Yamana-Residue-Periods-JFA}. The idea there is to apply the relative truncation operator only to $\varphi$. Following \cite{Ichino-Yamana-Compositio-Periods,Yamana-Residue-Periods-JFA} and work of Zydor \cite{Zydor-periods-AENS}, we can regularize the integral \eqref{formula::intro--Rankin-selberg-integral}. A key property that we need is the analogue of Theorem 3.10 in \cite{J-L-R-periods-JAMS} that expresses the period integral of truncated automorphic forms in terms of the regularized period integral with repect to a Levi subgroup of the constant terms. We will eventually obtain an identity from this property of the form
	\begin{align}\label{formula::intro--key property}
		\int^{\ast}_{G'(F)\backslash G'(\BA)} \varphi(g) dg   =  \int^{\ast}_{G'(F)\backslash G'(\BA)} \varphi(g) E(g,\Phi,s) dg  + \left(\sum_P \int^{\ast}_{M'(F)\backslash M'(\BA)} \cdots \right),
	\end{align}
	where the superscript $\ast$ denotes for the regularization of the integral and the remainder terms are taken over certain proper parabolic subgroups. Once one can deal with the first term in the right hand of \eqref{formula::intro--key property} in some other way, one can conclude information for the regularized integral of $\varphi$ by induction. 
	
	We switch to our specific situation. Let $G = G_{2n} = \GL_{2n}$ and $G' = G'_{2n} \cong \GL_n \times \GL_n$ be the centralizer of the element
	\begin{align*}
		\epsilon_{2n} = \left(\begin{smallmatrix}
			1 &  & & & &  \\
			& -1 & & & &  \\
			& & 1 & & & \\
			& & & -1 & & \\
			& & & & \ddots & \\
			& & & & & -1
		\end{smallmatrix} \right).
	\end{align*}
	The periods in this case are the so-called linear periods first studied by Friedberg and Jacquet \cite{Jacquet-Friedberg-Crell-LinearPeriods} (also called Friedberg-Jacquet periods in literature). It is closely related to the Bump-Friedberg $L$-function which is a product of the standard $L$-function and the exterior square $L$-function, see \cite{Bump-Friedberg-Exterior-L,Jacquet-Friedberg-Crell-LinearPeriods,Matringe-specialisation-BF-L}. Regularization is indispensible for the definition of linear periods of square-integrable automorphic forms on $G(\BA)$ as seen evidently from constant functions. However, we remark that not all linear periods of square-integrable automorphic forms can be regularized in the sense of Zydor. To remedy this, we take an analytic family of characters of $G'(\BA)$ and consider the linear periods with respect to these characters. The twisted period is then a holomorphic function of the analytic parameter whenever it is well defined. We then prove that it can be extended to a holomorphic function on all the complex plane.

	The discrete spectrum of $G(\BA)$ was classified by Moeglin and Waldspurger in \cite{Moeglin-Waldspurger-residue-spectrum-GLn}. For a divisor $r$ of $2n$ with $rd = 2n$ and an irreducible cuspidal automorphic representation $\sigma$ of $G_r(\BA)$, let $\speh(\sigma,d)$ be the representation in the discrete spectrum of $G(\BA)$ corresponding to $(\sigma,r)$, as well as its realization in $L^2(G(F)\backslash G(\BA))$ as iterated residues of Eisenstein series. 
	
	Let $G'(\BA)^{1,G} = \{ g \in G'(\BA) \ |\ |\det g| =1\}$. For $s \in \BC$, let $\mu_s$ be the character on $G'(\BA)$ defined by
	\begin{align*}
		\mu_s (\iota(g_1,g_2)) = |\det g_1|^s |\det g_2|^{-s},
	\end{align*} 
	where $\iota(g_1,g_2) \in G'(\BA)$ is given in \eqref{formula::defn-element-G'}. Let $\psi$ be a fixed nontrivial additive character of $F\backslash \BA$, also viewed as a non-degenerate character of the maximal unipotent subgroup of $G$ in the ususal way. We refer to Section \ref{sec::cusp linear periods-notation} for unexplained notation in the following theorem.
	
	\begin{thm}\label{thm::1}
		Let $2n = dr$. Let $\sigma$ be an irreducible cuspidal representation of $G_r(\BA)$ with trivial central character. Let $Q$ denote the standard parabolic subgroup of $G$ of type $(2n-r,r)$. If $r$ is odd, then for all $\varphi \in \speh(\sigma,d)$,
		\begin{align*}
			\int^{\ast}_{G'(F)\backslash G'(\BA)^{1,G}} \varphi(g) \mu_{s}(g) dg = 0.
		\end{align*}		
		If $r$ is even and $d \geqslant 2$, then for all $\varphi \in \speh(\sigma,d)$, 
		\begin{align}\label{formula::Theorem 1--I}
			\int^{\ast}_{G'(F)\backslash G'(\BA)^{1,G}}  \varphi(g) \mu_s (g) dg
		\end{align}
		is equal to 
		\begin{align}\label{formula::Theorem 1--II}
		rv_Q	\int_{P_n^{(r/2)}(\BA) \times P_{n-1}^{(r/2-1)}(\BA) \backslash G_n(\BA) \times G_{n-1}(\BA)} \int_{G'_{2n-r}(F) \backslash G'_{2n-r}(\BA)^{1,G_{2n-r}}}^{\ast}    \int_{N_r(F) \backslash N_r(\BA)}  \\
			 \varphi_Q \left( \begin{pmatrix}
				m  &  \\  & u 
			\end{pmatrix}  g \right) \overline{\psi(u)} \mu_{s-1/2}(m)  \mu_s(g) du dm dg. \nonumber
		\end{align}
		Here $v_Q$ is a certain volume (see \eqref{formula::defn-v_P-volume}), $G_n(\BA) \times G_{n-1}(\BA)$ is viewed as a subgroup of $G'_{2n}(\BA)$ and $\varphi_Q$ is the constant term along $Q$.
	\end{thm}
	The measure on the quotient $P_n^{(r/2)}(\BA) \times P_{n-1}^{(r/2-1)}(\BA) \backslash G_n(\BA) \times G_{n-1}(\BA)$ is right invariant. The integral in \eqref{formula::Theorem 1--II} formally makes sense (see Lemma \ref{lem::constant term--residue}).
	The equality of \eqref{formula::Theorem 1--I} and \eqref{formula::Theorem 1--II} holds in the following sense. \eqref{formula::Theorem 1--II} has a factorization into a product of local integrals. The local integrals and the infinite product converge absolutely and the product of local integrals is also a factorization of \eqref{formula::Theorem 1--I}. 
	
	 Our formula for the periods of sqaure-integralbe automorphic forms has an inductive feature. This shows by induction that the periods \eqref{formula::Theorem 1--I} is holomorphic at all $s \in \BC$. So we make the following definition.
	 
	 \begin{defn}
	 	An irreducible, discrete representation $\pi$ of $G(\BA)$ is called $(G',\mu_s)$-distinguished if there is an automorphic form $\varphi$ in the space of $\pi$ such that the period integral
	 	\begin{align*}
	 		\int^{\ast}_{G'(F)\backslash G'(\BA)^{1,G}}  \varphi(g) \mu_s (g) dg
	 	\end{align*}
 		is nonzero.
	 \end{defn}
	 The last equality in Theorem \ref{thm::1} realizes the regularized twisted period on $G_{2n-r}(\BA)$ as an inner integral of the regularized twisted period on $G_{2n}(\BA)$. By induction and an argument in the local theory of integral representations, we prove the following theorem.
	 
	\begin{thm}\label{thm::2}
		Let $\sigma$ be an irreducible cuspidal automorphic representation of $G_r(\BA)$ with trivial central character. Then $\speh(\sigma,d)$ is $(G',\mu_s)$-distinguished if and only if $r$ is an even integer and $\sigma$ is $(G_r',\mu_s)$-distinguished resp. $(G'_r,\mu_{s-1/2})$-distinguished when $d$ is even resp. $d$ is odd.
	\end{thm}
	
	In this paper we also discuss the case $(\GL_n,\GL_p \times \GL_q)$ where $p \ne q$. As a direct consequence of the study of local periods, the global regularized twisted periods vanishes identically.
	\begin{thm}
		Let $n = dr$ and $\sigma$ be an irreducible cuspidal automorphic representation of $G_r(\BA)$ with trivial central character. Let  $\varphi \in \speh(\sigma,d)$. Then
		\begin{align*}
			\int^{\ast}_{[G']^{1,G}}\varphi (g) \xi_s (g) dg =  0.
		\end{align*}
	\end{thm}

	The paper is organized as follows. In Section \ref{sec::notation} we set up notation. In Section \ref{sec::regularization} we develop the regularization of periods of automorphic forms on $G \times G'$ when $G'$ is a reductive subgroup of $G$, following Ichino and Yamana's work \cite{Ichino-Yamana-Compositio-Periods} and Zydor's work \cite{Zydor-periods-AENS}. We feel no need to confine ourselves to our specific case. Most statements and proofs in this section are nearly identical with those in \cite{J-L-R-periods-JAMS,Lapid-Rogawski-periods-Galois,Ichino-Yamana-Compositio-Periods,Zydor-periods-AENS} except that we make an extra assumption on the pair $(G,G')$ when considering the periods of truncated automrophic forms in Section \ref{sec::truncated periods} (see Remark \ref{rmk::Zydor's work}). In Section \ref{sec::cusp linear periods-notation} we recall the Rankin-Selberg theory of linear periods for cuspidal automorphic forms. In Section \ref{sec::local periods} we focus on the local theory which are key to the vanishing results. In Section \ref{sec::Eisenstein series-const terms} and Section \ref{sec::linear period-square-integral case} we recall some preliminaries on the theory of Eisenstein series and study the regularized periods of square-integrable automorphic forms. We prove Theorem 1.1 there. We prove Theorem 1.2 in Section \ref{sec::dist-spectrum} and discuss the case $p\ne q $ in Section \ref{sec:: case-p not q}.

	\section{Notation and preliminaries}\label{sec::notation}
	
	Let $F$ be a number field with adele ring $\BA$. Let $G$ be a connected reductive algebraic group over $F$. All subgroups are assumed to be closed and defined over $F$. Fix $A_0$ a maximal $F$-split torus of $G$, its centralizers $M_0$ in $G$ and $P_0$ a minimal parabolic group containing $M_0$. The standard resp. semi-standard parabolic subgroups are those containing $P_0$ resp. $A_0$. Denote by $\CF(A_0)$ the subset of all semi-standard parabolic subgroups of $G$.
	
	For a semi-standard parabolic subgroup $P$, when writing $P = MU$, we mean $U = U_P$ the unipotent radical of $P$ and $M = M_P$ the Levi subgroup of $P$ containing $A_0$. Denote by $A_P$ the split center of $M=M_P$. It is the subtorus of $A_0$ centralizing $M$. For an algebraic group over $F$, let $\rat(H)$ be the group of algebraic characters of $H$ defined over $F$. Set
	\begin{align*}
		\mathfrak{a}_P^{\ast} = \mathfrak{a}_M^{\ast} = \rat(M) \otimes_{\BZ} \BR, \quad \mathfrak{a}_P = \mathfrak{a}_M = \Hom_{\BZ}(\rat(M),\BR).
	\end{align*}
	We sometimes use $P$ or $M$ interchangeably in subscripts or superscripts. We write $\mathfrak{a}_0$ for $\mathfrak{a}_{P_0}$. We have $\rat(M)\otimes_{\BZ} \BR \cong \rat(A_M) \otimes_{\BZ} \BR$. For semi-standard parabolic subgroups $P \subset Q$, we have naturally $\mathfrak{a}_Q \subset \mathfrak{a}_P \subset \mathfrak{a}_0$. Fix a $W$-invariant inner product $\langle \cdot,\cdot \rangle$ on $\mathfrak{a}_0$, where $W$ is the Weyl group associated to $(G,A_0)$. We can then identify $\mathfrak{a}_P$ with $\mathfrak{a}_P^{\ast}$. 
	
	The inner product induces a topology on $\mathfrak{a}_0$. We write $\bar{A}$ for the closure of a subset $A $ of $\mathfrak{a}_0$ in this topology.
	 
	For semi-standard parabolic subgroups $P \subset Q$, let $\mathfrak{a}_P^Q$ be the orthogonal complement of $\mathfrak{a}_Q$ in $\mathfrak{a}_P$ and $\Delta_P \subset \mathfrak{a}_P$ be the set set of simple roots for the action of $A_P$ on $U_P$. Let $\Delta_P^Q$ be the subset of $\Delta_P$ that vanishes on $\mathfrak{a}_Q$. It is a basis of $\mathfrak{a}_P^Q$. Denote by $(\hat{\Delta}^{\vee})_P^Q$ the basis of $\mathfrak{a}_P^Q$ dual to $\Delta_P^Q$. We omit the superscript if $Q  = G$. Define
	\begin{align*}
		\mathfrak{a}_P^+  =  \{H \in \mathfrak{a}_P \ |\ \langle H,\alpha \rangle > 0 \text{ for all }\alpha \in \Delta_P \}.
	\end{align*}
 	
	Fix a good maximal compact subgroup $K$ of $G(\BA)$ that is adapted to $M_0$. Let $P \in \CF(A_0)$, define the Harish-Chandra function $H_P : M(\BA) \ra \mathfrak{a}_P$ by
	\begin{align*}
		e^{ \langle \chi,H_P(m) \rangle } = |\chi(m)|
	\end{align*}
	for all $\chi \in \rat(M)$ and $m \in M(\BA)$. We then extend $H_P$ to a left $U(\BA)$, right $K$-invariant function on $G(\BA)$ using the Iwasawa decomposition with respect to $P$. Let $M(\BA)^1$ be the insection of the kernels of the homomorphism $|\chi|$, where $\chi$ is taken over $\rat(M)$. Choose an isomorphism $A_P \cong \BG_m^l$, where $\BG_m$ is the multiplicative group. Define $A_P^{\infty}$ to be the image of $(\BR_+^{\times})^l$ in $A_{P,\infty}$, where $A_{P,\infty}$ is the archimedean component of $A_P(\BA)$, $l = \dim \mathfrak{a}_P$ and $\BR \ra F \otimes_Q\BR$ is given by $x \mapsto 1 \otimes x$. The map $H_P$ induces an isomorphism $A_P^{\infty} \cong \mathfrak{a}_P$. We denote by $e^X$ the element in $A_P^{\infty}$ such that $H_P(e^X) = X$. Let $\rho_P$ be half the sum of characters $A_P$ acting on $N_P$. Hence the modular character $\delta_{P(\BA)}$ on $P(\BA)$ is given by $e^{\langle 2 \rho_P, H_P(\cdot)\rangle}$.	
	
	Let $\CA_P(G)$ be the space of automorphic forms on $U(\BA)M(F) \backslash G(\BA)$, see \cite[Section I.2.17]{Moeglin-Waldspurger-95-spectral}. For any locally integrable function $\phi$ on $P(F) \backslash G(\BA)$ and any parabolic subgroup $Q  \subset P$, the constant term of $\phi$ along $Q$ is defined by 
	\begin{align*}
		\phi_Q (g)  =  \int_{V(F) \backslash V(\BA)} \phi (vg) dv,
	\end{align*}
	where $V$ is the unipotent radical of $Q$. According to \cite[Section I.3.2]{Moeglin-Waldspurger-95-spectral}, every automorphic form $ \phi \in \mathcal{A}_P (G)$ has a finite sum decomposition
	\begin{align}\label{formula::notion--decomposition-autoForm}
		\phi ( g )  =  \sum_i  q_i (H_P(g))  e^{ \langle \lambda_i + \rho_P, H_P(g)\rangle} \phi_i (g)
	\end{align} 
	where $q_i \in \BC[\mathfrak{a}_P]$, $ 0 \neq \lambda_i \in \mathfrak{a}_{P,\BC}^*$ and $ 0 \neq \phi_i \in \mathcal{A}_P(G)$ such that $\phi_i(ag) = \phi_i (g)$ for $a \in A_P^{\infty}$ and $g \in G(\BA)$. The set composed of distinct $\lambda_i$ is uniquely determined by $\phi$ and is called the set of exponents of $\phi$. For $Q \subset P$ the set of exponents resp. cuspidal exponents of $\phi$ along $Q$ is by definition the set of exponents of $\phi_Q$ resp. $\phi^{\textup{cusp}}_Q$ and is denoted by $\mathcal{E}_Q(\phi)$ resp. $\mathcal{E}_Q^{\textup{cusp}}(\phi)$.

	\section{Periods of automorphic forms on $G \times G'$}\label{sec::regularization}
	
	Let $G' \subset G$ be a connected reductive subgroup. In \cite{Zydor-periods-AENS}, Zydor have defined the regularized period along $G'$ of an automorphic form on $G$ under some restrictions on its exponents, extending poineer works \cite{J-L-R-periods-JAMS} and \cite{Lapid-Rogawski-periods-Galois}. In this section we adapt this regularization procedure to the group $G \times G'$ and its diagonal subgroup $G'$, following the pattern laid down for Rankin-Selberg case and Galois symmetric case in \cite{Ichino-Yamana-Compositio-Periods} and \cite{Yamana-Residue-Periods-JFA}.
	
	\subsection{Relative truncation operator of Zydor} 
	
	We recall some notations from \cite{Zydor-periods-AENS}. Fix $A_0'$ a maximal $F$-split torus of $G'$ and $A_0$ a maximal $F$-split torus of $G$ such that $A_0' \subset A_0$. Also fix $P_0'$ a  minimal parabolic subgroup of $G'$ that contains $A_0'$. Denote by $\CF^G(A_0)$ resp. $\CF^{G'}(A_0')$ the set of parabolic subgroups of $G$ resp. $G'$ containing $A_0$ resp. $A_0'$. The family of parabolic subgroups that is pertinent to the regularizaton procedure is 
	\begin{align}\label{formula::Zydor--relavent parabolics}
		\CF^{G}(P_0') = \{ P \in \CF^G(A_0) \ |\ \overline{\mathfrak{a}_{0'}^+} \cap \mathfrak{a}_P^+ \ne \emptyset \}.
	\end{align}
	For all $P \in \CF^G(P_0')$, we will write $P' = P \cap G'$. Then $P'$ is a standard parabolic subgroup of $G'$ (containing $P_0'$), and 
	\begin{align}
		\mathfrak{a}_{P'}^+  \cap \mathfrak{a}_P^+ \ne \emptyset,
	\end{align}
	see \cite[Proposition 3.1]{Zydor-periods-AENS}. Let $P = M \ltimes U$ and $P' = M' \ltimes U'$ be the Levi decomposition of $P$ and $P'$ respectively. Then
	\begin{align*}
		M' = M \cap G',\quad U' = U \cap G'.
	\end{align*}
	Using dynamic description of parabolic subgroups, $P$ can be defined by a cocharacter $\mathbb{G}_m \ra G' \subset G$ that defines $P'$ in $G'$.
	
	Following \cite[Section 3.2]{Zydor-periods-AENS}, we set for $P \in \CF^G(P_0')$ 
	\begin{align*}
		\mathfrak{z}_P = \mathfrak{a}_P \cap \mathfrak{a}_{0'}, \quad \mathfrak{z}_{P}^+ = \mathfrak{a}_P^+ \cap \mathfrak{a}_{0'}.
	\end{align*}
	We can and will assume that the Euclidean, Weyl group invariant structure on $\mathfrak{a}_{0'}$ is induced from the one on $\mathfrak{a}_0$. For $P,\,Q \in \CF^G(P_0')$ such that $P \subset Q$, clearly we have $\overline{\mathfrak{z}_Q^+} \subset \overline{\mathfrak{z}_P^+}$. We denote by $\mathfrak{z}_P^Q$ the orthogonal complement of $\mathfrak{z}_Q$ in $\mathfrak{z}_P$, and let 
	$$\varepsilon_P^Q = (-1)^{\dim \mathfrak{z}_P^Q}.$$
	Denote $\mathfrak{z}^P$ the orthogonal complement of $\mathfrak{z}_P$ in $\mathfrak{a}_{0'}$. For $X \in \mathfrak{a}_{0'}$, denote by $X_P$, $X^P$ and $X_P^Q$ projections of $X$ onto $\mathfrak{z}_P$, $\mathfrak{z}^P$ and $\mathfrak{z}_P^Q$ respectively.
	
	For $P,Q \in \CF^G(P_0')$ such that $P \subset Q$, take $z \in \mathfrak{z}_Q^+$. As in \cite[Section 3.2]{Zydor-periods-AENS}, let $\tau_P^Q$ be the characteristic function of the relative interior of the cone
	\begin{align*}
		A(\overline{\mathfrak{z}_Q^+},\overline{\mathfrak{z}_P^+}) : =\{ \lambda ( x - z) \ |\ \lambda >0, x \in \overline{\mathfrak{z}_P^+} \} \subset \mathfrak{z}_P,
	\end{align*} 
	and let $\hat{\tau}_P^Q$ be the characteristic function of the relative interior of 
	\begin{align*}
		A(\overline{\mathfrak{z}_Q^+},\overline{\mathfrak{z}_P^+})^{\vee} : = \{ X \in \mathfrak{a}_{0'} \ |\ \langle X, x -z \rangle \geqslant 0 \text{ for all }x\in \overline{\mathfrak{z}_P^+} \} \subset \mathfrak{z}^Q.
	\end{align*}
	The two definitions are independant of the choice of $z$. We omit the superscript if $Q = G$. Observe that, for $X \in \mathfrak{a}_{0'}$, 
	\begin{align*}
		\tau_P^Q (X )  = \tau_P^Q (X^Q), \quad \hat{\tau}_P^Q (X) = \hat{\tau}_P^Q (X_P).
	\end{align*}
	A generalization of the Langlands combinatorial lemma \cite[Proposition 1.5]{Zydor-periods-AENS} asserts that for any $P \in \CF^G(P_0')$ and $X \in \mathfrak{a}_{0'}$ we have
	\begin{align}\label{formula::combinatorial lemma}
		\sum_{\stackrel{Q \supset P}{P,Q\in \CF^G(P_0')}}  \varepsilon_P^Q \hat{\tau}_P^Q (X^Q) \tau_Q (X_Q)    = \begin{cases*}
			1 \quad \text{if }P=G \\
			0 \quad \text{otherwise}. 
		\end{cases*}
	\end{align} 
	
	For $P,Q \in \CF^G(P_{0'})$ such that $P \subset Q$, as in \cite[Section 3.3]{Zydor-periods-AENS}, we let
	\begin{align}\label{formula::defn-Gamma}
		\Gamma_P^Q (H,X) = \sum_{\stackrel{R \in \CF^G(P'  _{0})}{P \subset R \subset Q}}  \varepsilon_R^Q \tau_P^R(H) \hat{\tau}_R^Q (H - X), \quad H,X \in \mathfrak{a}_{0'}.
	\end{align}
	For all $X \in \mathfrak{a}_{0'}$, the function $\Gamma_P^Q(\cdot,X)$ is compactly supported \cite[Lemma 3.3]{Zydor-periods-AENS}.
	
	Let $Q \in \mathcal{F}^G(P'_0)$ and $Q' = Q \cap G'$. Let $\phi$ be a locally integrable function on $Q(F)\backslash G(\BA)$. The relative truncation operator of Zydor is defined as follows \cite[Section 3.7]{Zydor-periods-AENS}
	\begin{align*}
		\Lambda^{T,Q}\phi(x) = \sum_{\stackrel{P\subset Q}{P\in \mathcal{F}^G(P_0')}} \varepsilon^Q_P \sum_{ \delta \in P'(F)\backslash Q'(F)} \hat{\tau}_P^Q (H_{0'}(\delta x)^Q - T^Q) \phi_P (\delta x),
	\end{align*}
	for all $x \in Q'(F) \backslash G'(\BA)$. When $Q = G$ we write $\Lambda^T$ for $\Lambda^{T,G}$. The sums in the definition of $\Lambda^{T,Q}$ are finite \cite[Lemma 2.8]{Zydor-periods-AENS}. There is a $\mathcal{T}_{reg} \in \mathfrak{a}_{0'}$ defined using the reduction theory for $G'$ \cite[Section 3.5]{Zydor-periods-AENS}. We will call elements in $\mathcal{T}_{reg} + \mathfrak{a}_{0'}$ sufficiently positive. For sufficiently positive $T \in \mathfrak{a}_{0'}$, the operator $\Lambda^{T,Q}$ carries smooth functions on $Q(F) \backslash G(\BA)$ of uniform moderate growth to functions on $V'(\BA)L'(F) \backslash G'(\BA)$ of rapidly decay \cite[Theorem 3.9]{Zydor-periods-AENS}. Here $Q' = L'V'$ is the Levi decomposition of $Q'$. If $\phi$ is a locally integrable function on $G(F) \backslash G(\BA)$, then clearly
	\begin{align*}
		\Lambda^{T,Q} \phi = \Lambda^{T,Q} \phi_Q.
	\end{align*}
	
	We have the inversion formula 
	\begin{align*}
		\phi(x) = \sum_{P\in \mathcal{F}^G(P_0')} \sum_{\delta \in P'(F)\backslash G'(F)} \Lambda^{T,P}\phi(\delta x) \tau_P(H_{0'}(\delta x)_P - T_P)
	\end{align*}
	for any locally integrable function $\phi$ on $G(F)\backslash G(\BA)$ and $x \in G'(\BA)$. This follows directly from \eqref{formula::combinatorial lemma}. Also, we have 
	\begin{align}\label{formula::lambda-T-T'}
		\Lambda^{T+T',P}\phi(x) = \sum_{\stackrel{Q\in \mathcal{F}^G(P_0')}{Q\subset P}} \sum_{ \delta \in Q'(F)\backslash P'(F)} \Lambda^{T,Q}\phi(\delta x)\Gamma^P_Q(H_{0'}(\delta x)^P_Q - T^P_Q, T'^P_Q ),
	\end{align}
	for any locally integrable function $\phi$ on $P(F) \backslash G(\BA)$, all $T,T' \in \mathfrak{a}_{0'}$ and $x \in G'(\BA)$.
	
	 Unless otherwise specified, when doing summation over parabolic subgroups, we shall always assume the parabolic subgroups are contained in $\CF^G(P_0')$. 
	
	Let $P \in \CF^G(P_{0'})$ and $P' = P \cap G'$. Let $P' = M'U'$ be the Levi decomposition of $P'$ with $M'$ containing $A_0'$. Set as in \cite[Section 3.4]{Zydor-periods-AENS}
	\begin{align*}
		M'(\BA)^{1,P} = \{ x \in M'(\BA) \ |\ H_{P'}(x)_P = 0  \}.
	\end{align*}
	Set $Z_P^{\infty} = A_P^{\infty} \cap A_{P'}^{\infty}$. The restriction of $H_{P'}$ to $Z_P^{\infty}$ is a group isomorphism with its image $\mathfrak{z}_P$. We then have a decomposition 
	\begin{align*}
		M(\BA) = Z_P^{\infty} M(\BA)^{1,P}.
	\end{align*}
	We will write $[G']^{1,G}$ for the relative automorphic quotient $G'(F) \backslash G'(\BA)^{1,G}$.
	
	
	\subsection{Integrals over cones}
	
	Let $V$ be a finite dimensional Euclidean space over $\BR$. Let $\langle \cdot,\cdot \rangle$ be the scalar product on it. Denote $V_{\BC} = V \otimes_{\BR} \BC$ and extend the bilinear product $\langle \cdot ,\cdot \rangle$ to it in a natural way. For $\lambda \in V_{\BC}$ we denote by $\Re \lambda$ the real part of $\lambda$.
	
	A polynomial exponential function on $V$ is a  function $f$ of the form
	\begin{align*}
		f(T) = \sum_{\lambda \in V_{\BC}} e^{ \langle \lambda ,T\rangle} q_{\lambda}(T)
	\end{align*}
	where $q_{\lambda} \in \BC[V]$ with $q_{\lambda} = 0$ for all but finitely many $\lambda \in V_{\BC}$. Such a decomposition is unique. The $\lambda$ such that $q_{\lambda} \ne 0$ are called the exponents of $f$. The purely polynomial part of $f$ is by definition the polynomial $q_0$ corresponding to $\lambda = 0$.
	
	By a cone in $V$ we shall mean a finite intersection of half-spaces of $V$. Here a half-space of $V$ is a subset of the form $\{X \in V \ |\ \langle X, H \rangle \geqslant 0\}$, for some nonzero $H$ in $V$. We refer the readers to \cite[Section 1]{Zydor-periods-AENS} and references therein for more about the theory of cones. By a face of $C$ we mean an intersection of $C$ with a half-space containing $C$ or $-C$. The set of faces of $C$ has a natural poset structure induced by inclusion relation. In particular, the minimal face of $C$ is the maximal subspace that is contained in $C$, and is denoted by $F_0(C)$. For a cone $C$, the dual cone $C^{\vee}$ is defined as
	\begin{align*}
		C^{\vee} = \{ X \in V \ |\ \langle X, C\rangle \geqslant 0 \}.
	\end{align*}
	We have $(C^{\vee})^{\vee} = C$. A cone $C$ in $V$ is called non-degenerate if it has a non-empty interior in $V$ and if $F_0(C)= \{0\}$. If $C$ is non-degenerate, then $C^{\vee}$ is also non-degenerate. Let $C$ be a non-degenerate cone and $\{R_1,\cdots,R_d\}$ be the set of $1$-dimensional faces of $C$ ($d \geqslant \dim V$). Choose $0 \ne v_i \in R_i$ for $1 \leqslant i \leqslant d$. Then
	\begin{align*}
		C = \BR_+ v_1 + \BR_+ v_2 + \cdots + \BR_+ v_d.
	\end{align*}
	For $\lambda \in V_{\BC}$, we say $\lambda$ is non-degenerate resp. negative with respect to $C$ if $\langle \Re \lambda ,v_i \rangle \ne 0$ resp. $<0$ for $i = 1,\cdots,d$. For an arbitrary cone $C$ in $V$, let $V_C$ be the subspace of $V$ generated by $C$ and $V_C^{F_0}$ the orthogonal complement of $F_0(C)$ in $V_C$. Then $V_C^{F_0} \cap C$ is a non-degenerate cone in $V_C^{F_0}$.
	
	The $\sharp$-integral of polynomial exponential functions over \textit{simplicial} cones ($d = \dim V$) has been discussed in \cite{J-L-R-periods-JAMS}. Many discussions there work verbatimly to the case of non-degenerate cones. For a set $S$, denote by $\tau^S$ the characteristic function of $S$. Suppose that $C$ is a non-degenerate cone in $V$. Let $f$ be a polynomial exponential function on $V$. The $\sharp$-integral
	\begin{align*}
		\int_V^{\sharp} f(x ) \tau^C(x-T) dx
	\end{align*} 
	exists if the exponents of $f$ are all non-degenerate with respect to $C$. We also say that $f(x)\tau^C(x-T)$ is $\sharp$-integrable. Note that the necessarity here may not hold for non-simplicial $C$. Also, the function
	\begin{align*}
		T \mapsto \int_V^{\sharp} f(x) \tau^C(x-T) dx
	\end{align*}
	is a polynomial exponential function in $T$ with the same exponents as $f$.
	
	Let $V = W_1 \oplus W_2$ be an orthogonal decomposition of $V$. Let $C_1$ and $C_2$ be non-degenerate cones of $W_1$ and $W_2$. Then $C = C_1 \oplus C_2$ is non-degenerate in $V$. Write $T = T_1 + T_2$ and $x = w_1 + w_2$ relative to $V = W_1 \oplus W_2$. If the exponents of $f$ are non-degenerate with respect to $C$, then
	\begin{align}\label{formula::combine-integrals-I}
		\int_V^{\sharp} f(x) \tau^C(x &- T ) dT \\
		&= \int_{W_2}^{\sharp} \left( \int_{W_1}^{\sharp} f(w_1 + w_2)\tau^{C_1}(w_1 - T_1) dw_1  \right) \tau^{C_2}(w_2 - T_2) dw_2. \nonumber
	\end{align}
	Let $V  =  W_1 \oplus W_2$ and $C_2$ be as above. Let $g(x)$ be a compactly supported function on $W_1$. Write $T = T_1 + T_2$ and $x = w_1 + w_2$ as above. Functions on $V$ of the form 
	\begin{align*}
		g(w_1 - T_1)\tau^{C_2}(w_2 - T_2)
	\end{align*}
	are called functions of type (C). If the restrictions of the exponents of $f$ to $W_2$ is non-degenerate with respect to $C_2$, then
	\begin{align}\label{formula::combine-integrals-II}
		\int_V^{\sharp} f(x) & g(w_1 -T_1) \tau^{C_2}(w_2-T_2) dx \\
		&= \int_{W_2}^{\sharp} \left( \int_{W_1} f(w_1 + w_2) g(w_1 - T_1) dw_1\right) \tau^{C_2}(w_2 - T_2) dw_2. \nonumber
	\end{align}

	For $i = 1,\cdots,r$ let $V = W_{i1} \oplus W_{i2}$ be an orthogonal decomposition of $V$ and let $C_{i2}$ be a non-degenerate cone in $W_{i2}$. Set $G_i(w_1 + w_2) = g_i(w_1) \tau^{C_{i2}}(w_2) $ for $w_j \in W_{ij}$, where $g_i$ is either a compactly supported function on $W_{i1}$ for each $i$ or a characteristic functon $\tau^{C_{i1}}$ with $C_{i1}$ a non-degenerate cone in $W_{i1}$ for each $i$.
	\begin{lem}\label{lem::Lem 6}
		Let $C$ and $C^{\ast}$ be non-degenerate cones in $V$. Let $f$ be a polynomial exponential function on $V$. Assume that $C_{i2},C \subset C^{\ast}$ for all $i$ or $C,C_{i1} \times C_{i2} \subset C^{\ast}$ for all $i$, depending on $g_i$. Assume that 
		\begin{align*}
			\tau^C(x) = \sum_{i=1}^r a_i G_i(x)
		\end{align*}
		for some constants $a_i$. Assume that $f(x)G_i(x - T)$ is $\sharp$-integrable for each $i$. Then $f(x) \tau^C(x- T)$ is $\sharp$-integrable and 
		\begin{align*}
			\int_V^{\sharp} f(x) \tau^C(x - T) dx = \sum a_i \int_V^{\sharp} f(x) G_i(x - T) dx.
		\end{align*}
	\end{lem}
	
	\subsection{Regularization of periods}\label{sec::regularization of periods}
	
	Before proceding to define the regularized integral, we present a lemma which reveals a descent structure on the set $\CF^G(P_0')$ of relative parabolic subgroups. 
	\begin{lem}\label{lem::rel parabolics-descent}
		Let $P = MU \in \CF^G(P_0')$. The map
		\begin{align*}
			\{ Q \subset P \ |\ Q \in \mathcal{F}^G(P'_0) \} & \ra \mathcal{F}^M( M' \cap P'_0)   \\
			Q  & \mapsto Q \cap M
		\end{align*}
		is a bijection.
	\end{lem}
	\begin{proof} 
		We first show that for $Q \in \mathcal{F}^G(P'_0)$, $Q \cap M  \in \mathcal{F}^M(M' \cap P'_0)$. By definition, there exists a cocharacter $\lambda \in X_*(A'_0)$ such that $Q = P_G (\lambda)$ and $Q \cap G' = P_{G'}(\lambda)$ is a standard parabolic subgroup of $G'$. Hence $P_{M'}(\lambda)$ is a standard parabolic subgroup of $M'$ and $Q \cap M = P_M(\lambda) \in \mathcal{F}^M(M' \cap P'_0)$. The map is injective as we have $Q = (Q \cap M) U$ when $Q \subset P$. 
		
		Given $R \in \mathcal{F}^M(M'\cap P'_0)$. By definition, there exists $\lambda_1 \in X_*(A'_0)$ such that $R = P_M (\lambda_1)$ and $R \cap M' = P_{M'}(\lambda_1)$ is a standard parabolic subgroup of $M'$. Take a nonzero $\mu \in X_*(A'_0)$ such that $\mu \in \mathfrak{a}^+_P \cap \mathfrak{a}^+_{P'}$. The existence of such $\mu$ is shown in \cite[Proposition 3.1]{Zydor-periods-AENS}. Let $N$ be a sufficiently large positive integer. Then $P_M(\lambda_1 + N \mu) = R$. Let $Q = P_G(\lambda_1 + N \mu)$. Note that for any $\alpha \in \Phi(A_0,U)$, $\langle \alpha , \lambda_1 + N \mu\rangle > 0$. Hence $Q \subset P$. It remains to show that $P_{G'}(\lambda_1 + N \mu)$ is standard in $G'$. The assertion follows as $\Phi^+(A'_0,G') = \Phi(A'_0,U') \cup \Phi^+(A'_0,M')$.
	\end{proof}
	
	For $P \in \CF^G(P_{0'})$ and $P' = P \cap G'$, we have the Iwasawa decomposition of $G'(\BA)^{1,G}$ relative to $P'$:
	\begin{align}\label{formula::Iwasawa-P}
		G'(\BA)^{1,G} =  U'(\BA) Z_P^{\infty,G}  M'(\BA)^{1,P}K',
	\end{align}
	where $Z_P^{\infty,G}$ is the subgroup of $Z_P^{\infty}$ corresponding to $\mathfrak{z}_P^G$ under $H_{P'}$. For $X \in \mathfrak{z}_P^G$, we shall always write $e^X \in Z_P^{\infty,G}$ for the inverse image of $X$ under $H_{P'}$. Observe that $G'(\BA)^{1,G}$ is unimodular. The integration formula corresponding to the decomposition \eqref{formula::Iwasawa-P} is as follows:
	\begin{align}\label{formula::integral formula--before regularization}
		\begin{aligned}
			&\int_{P'(F) \backslash G'(\BA)^{1,G}} f(g) dg \\
			& = \int_{K'} \int_{M'(F) \backslash M'(\BA)^{1,P}} \int_{\mathfrak{z}_P^G} \int_{N'(F) \backslash N'(\BA)} f (n e^X m k )  e^{-\langle 2\rho_{P'}, H_{P'}(e^X m)\rangle} dn dx dm dk,
		\end{aligned}
	\end{align}
	for any $f$ locally integrable over $P'(F) \backslash G'(\BA)^{1,G}$.

	Let $P \in \CF^G(P_{0'})$ and $P' = P \cap G'$. Automorphic forms $\phi \in \CA_P(G)$ and $\phi' \in \CA_{P'}(G')$  have decompositions of type \eqref{formula::notion--decomposition-autoForm}. Write $g = u e^X m k$ for an Iwasawa decomposition of $g \in G'(\BA)^{1,G}$  with respect to $P'$, where $u \in U'(\BA)$, $ X \in \mathfrak{z}_P^G$, $ m \in M'(\BA)^{1,P}$ and $k \in K'$. If follows from \eqref{formula::notion--decomposition-autoForm} that we can write
	\begin{align}\label{formula::expansion--before regularization}
		\begin{aligned}
			&\phi( u e^X m k)  = \sum_i q_i (X)e^{\langle \lambda_i + \rho_P, X  \rangle} \phi_i (mk),\\
			&\phi'(u e^X m k ) = \sum_j  q'_j (X)e^{\langle \lambda_j'+ \rho_{P'}, X  \rangle } \phi'_j (mk).
		\end{aligned}
	\end{align}
	where 
	\begin{align*}
		\begin{aligned}
			&q_i  \in \BC[\mathfrak{z}_P^G],\  \lambda_i \in \CE(\phi),\ \phi_i \in \CA_P(G),\\
			&q'_j \in \BC[\mathfrak{z}_P^G],\ \lambda_j' \in \CE(\phi'),\ \phi_j' \in \CA_{P'}(G').
		\end{aligned}
	\end{align*}
	Note that, in general, $\phi_i$ and $\phi'_j$ are no longer left $Z_P^{\infty}$-invariant. It is not hard to check that
	\begin{align*}
		\Lambda^{T,P}\phi (g ) \phi' ( g ) = \sum_{i,j} & q_i ( X ) q'_j ( X ) e^{\langle \lambda_i + \rho_P + \lambda'_j + \rho_{P'} ,X \rangle}  \Lambda^{T,P} \phi_i(mk) \phi'_j(mk).
	\end{align*}
	
	Let $\tau_k(X)$ be a function of type (C) on $\mathfrak{z}_P^G$ that depends continuously on $k \in K'$ in the sense of \cite[Section 7]{J-L-R-periods-JAMS}.

	Recall that for $P \in \CF^G(P_{0'})$, $\underline{\rho}_P$ is the projection of $\rho_P - 2 \rho_{P'}$ onto $\mathfrak{a}_{0'}$. We define the $\sharp$-integral
	\begin{align}\label{formula::regularized period--I}
		\int^{\sharp}_{P'(F) \backslash G'(\BA)^{1,G}} \Lambda^{T,P}\phi(g) \phi'(g) \tau_k (H_{0'}(g)_P - T_P) dg
	\end{align}
	as
	\begin{align}\label{formula::regularized period--II}
		\sum_{i,j} \int_{K'} &\left( \int_{M'(F) \backslash M'(\BA)^{1,P}}  \Lambda^{T,P} \phi_i(mk) \phi'_j(mk) e^{ -\langle 2\rho_{P'},H_{P'}(m) \rangle } dm \right)   \\
		  &\quad \quad \qquad \qquad \left(\int_{\mathfrak{z}_P^G}^{\sharp} q_i(X)q'_j(X) e^{ \langle \lambda_i  + \lambda'_j + \rho_{P'} + \underline{\rho}_P , X \rangle} \tau_k (X -T_P)  dX  \right)  dk, \nonumber
	\end{align}
	whenever it makes sense. The inner integral over $m \in M'(F) \backslash M'(\BA)^{1,P}$ is absolutely convergent since $\Lambda^{T,P}\phi_i$ is of rapidly decay \cite[Theorem 3.9]{Zydor-periods-AENS} and $\phi'_j$ is of moderate growth as functions on $M'(\BA)^{1,P}$. Thus, \eqref{formula::regularized period--I} is defined if all inner $\sharp$-integrals over $\mathfrak{z}_P^G$ in \eqref{formula::regularized period--II} exist. 
		
	\begin{defn}
		Let $\mathcal{A}(G \times G')^*$ be the space of pairs $(\varphi,\varphi') \in \mathcal{A}(G) \oplus \mathcal{A}(G')$ which satisfy
		\begin{align*}
			\langle \lambda +  \lambda' + \rho_{P'} + \underline{\rho}_P,\mathfrak{z}_Q^G \rangle  \neq 0
		\end{align*}
		for any $P \subset Q$ in $\CF^G(P'_0)$ with $\dim \mathfrak{z}_Q^G = 1$ and 
		\begin{align*}
			\lambda \in \mathcal{E}_P (\varphi), \quad \lambda' \in \mathcal{E}_{P'}(\varphi').
		\end{align*}
		When $(\varphi,\varphi') \in \mathcal{A}(G \times G')^*$, all the $\sharp$-integrals
		\begin{align*}
			\mathbf{P}^{G',T}_P (\varphi \otimes \varphi') = \int^{\sharp}_{P'(F)\backslash G'(\BA)^{1,G}} \Lambda^{T,P}\varphi(g)\varphi'_{P'}(g)  \tau_P (H_{0'}(g)_P - T_P)dg
		\end{align*}
		exist and a regularized period $\mathbf{P}^{G'}(\varphi \otimes \varphi')$ is defined as the sum
		\begin{align*}
			\sum_{P\in \mathcal{F}^G(P_0')} \mathbf{P}^{G',T}_P(\varphi \otimes \varphi').
		\end{align*}
	\end{defn}
	
	\begin{prop}\label{prop::period-well defined}
		\textup{(i)} \ $\mathbf{P}^{G'}$ is well-defined and is independent of $T$.
		
		\textup{(ii)}\  $\period^{G'}$ defines an $G'(\BA_f)^{1,G}$-invariant linear functional on $\CA(G \times G')^{\ast}$.
	\end{prop}
	\begin{proof}
		The proof of these two statements follows almost word for word the proof of \cite[Theorem 9]{J-L-R-periods-JAMS}. We are content with ourselves to write down only the details of the proof of \textup{(i)} in this general setting.  Let $P \in \CF^G(P'_0)$ and $P' = P \cap G'$. Let $g = u e^X m k$ be an Iwasawa decomposition of $g \in G'(\BA)^{1,G}$  relative to $P'$ as above. Choose decompositions of $\varphi_P$ and $\varphi'_{P'}$ as in \eqref{formula::expansion--before regularization}. 
		We use \eqref{formula::lambda-T-T'} to write $\period^{G',T+T'}_P (\varphi \otimes \varphi')$ as the integral over $k \in K'$ and sum over $i,j$ of the product of 
		\begin{align}\label{formula::period-well-defined-I}
			\int_{\mathfrak{z}_P^G}^{\sharp}q_i(X)q'_j(X) e^{ \langle \lambda_i  + \lambda'_j + \rho_P + \rho_{P'} - 2\rho_{P'}  , X \rangle} \tau_P (X -T_P - T'_P)  dX
		\end{align}
		and
		\begin{align}\label{formula::period-well-defined-II}
			\int_{M'(F) \backslash M'(\BA)^{1,P}} \sum_{Q\subset P}  \sum_{ \delta \in Q'(F)\backslash P'(F)} \Gamma^P_Q(H_{0'}(\delta m k)^P_Q - T^P_Q, &T'^P_Q ) \Lambda^{T,Q} \phi_i(\delta mk)\\
			&\cdot \phi'_j(mk)  e^{ -\langle 2\rho_{P'},H_{P'}(m) \rangle }dm. \nonumber
		\end{align}
		Taking the sum over $Q$ outside the integral, \eqref{formula::period-well-defined-II} is equal to the sum over $Q \subset P$ of 
		\begin{align}\label{formula::period-well-defined-III}
			\int_{(M'\cap Q')(F) \backslash M'(\BA)^{1,P}} \Gamma^P_Q(H_{0'}( m k)^P_Q - T^P_Q, T'^P_Q ) \Lambda^{T,Q} \phi_i( mk) \phi'_j(mk) e^{ -\langle 2\rho_{P'},H_{P'}(m) \rangle }dm.
		\end{align}
		The equality will be justified if we can show, for each $Q$, the integral \eqref{formula::period-well-defined-III}   is absolutely convergent. We may replace $\phi_i$ and $\phi'_j$ by $(\phi_i)_Q$ and $(\phi'_j)_{Q'}$. This will be clear once we have justified the convergence of \eqref{formula::period-well-defined-III}. Let $m = v e^Y l k_1$
		be an Iwasawa decomposition of $m \in M'(\BA)^{1,M}$ with respect to $M' \cap Q'$ with $v \in (M' \cap V')(\BA)$, $Y \in \mathfrak{z}_Q^P$, $l \in L'(\BA)^{1,M\cap Q}$ and $k_1 \in K' \cap M'(\BA)$. We  expand $(\phi_i)_Q$ and $(\phi_j')_{Q'}$ as in \eqref{formula::expansion--before regularization}.
		\begin{align*}
			\begin{aligned}
				&(\phi_i)_Q ( v e^Y l k_1 k) = \sum_{s} q_{is} (Y)e^{\langle \lambda_{is} + \rho_Q,Y \rangle} \phi_{is} (l k_1 k)  \\
				&(\phi'_j)_{Q'} ( v e^Y l k_1 k)  = \sum_{t} q'_{jt} (Y) e^{\langle \lambda'_{jt}+ \rho_{Q'},Y\rangle} \phi'_{jt} (l k_1 k).
			\end{aligned}
		\end{align*}
		 In view of the integration formula \eqref{formula::integral formula--before regularization} adapted to the current case, we rewrite \eqref{formula::period-well-defined-III} as the integral over $ k_1 \in K' \cap M'(\BA)$ and sum over $s,t$ of the product of 
		 \begin{align}\label{formula::period-well-defined-IV}
 	  		\int_{\mathfrak{z}_Q^P} \Gamma_Q^P (Y - T_Q^P, T'^P_Q) q_{is} (Y)q'_{jt} (Y)  e^{\langle \lambda_{is} +\lambda'_{jt} + \rho_Q + \rho_{Q'} - 2 \rho_{ Q'} ,Y \rangle} dY
		 \end{align}
	 	and 
	 	\begin{align}\label{formula::period-well-defined-V}
	 		\int_{[L']^{1,Q}} \Lambda^{T,Q}\phi_{is}(l k_1 k) \phi'_{jt}(l k_1 k)   e^{- \langle 2 \rho_{Q'}, H_{L'}(l) \rangle} dl.
	 	\end{align}
 		Here we have used the relation
 		\begin{align}\label{formula::rho-relation}
 			\rho_{P'} + \rho_{M'_{Q'}}  = \rho_{Q'}.
 		\end{align}
		Since $\Gamma_Q^P(Y-T_Q^P,T'^P_Q)$ is a compactly supported function of $Y \in \mathfrak{z}_Q^P$, the integral \eqref{formula::period-well-defined-IV} is absolutely convergent. By Fubini's theorem, we get that \eqref{formula::period-well-defined-III} is absolutely convergent. Since we are integrating over $K'$, we may drop the integral over $K' \cap M'(\BA)$. Thus, $\period_P^{G',T+T'}(\varphi \otimes \varphi')$ equals to the integral over $k \in K'$, sum over $Q \subset P$, and sum over $i,j,s,t$ of the product of \eqref{formula::period-well-defined-V}, \eqref{formula::period-well-defined-I} and \eqref{formula::period-well-defined-IV}. According to \eqref{formula::combine-integrals-II}, we may combine the integrals over $\mathfrak{z}_Q^P$ and $\mathfrak{z}_P^G$ to a $\sharp$-integral over $\mathfrak{z}_Q^G$. It is not hard to check that $\period_P^{G',T+T'}(\varphi \otimes \varphi')$ is equal to the sum over $Q \subset P$ of
		\begin{align}\label{formula::period-well-defined-VI}
			\int^{\sharp}_{Q'(F) \backslash G'(\BA)^{1,G}}  \Lambda^{T,Q} \varphi (g) \varphi'_{Q'}(g) \tau_P (H_{0'}(g)_P^G - T_P - T'_P) \Gamma_Q^P (H_{0'}(g)_Q^P - T_Q^P, T'^P_Q) dg.
		\end{align}
		Summing over $P$, we see that $\period^{G',T+T'}(\varphi\otimes \varphi')$ is the sum over pairs $Q \subset P $ of \eqref{formula::period-well-defined-VI}. As $\mathfrak{z}_P^G \cap \overline{\mathfrak{z}_P^+}$ is contained in $\mathfrak{z}_Q^G \cap \overline{\mathfrak{z}_Q^+}$, we may apply Lemma \ref{lem::Lem 6} to take $P$ inside the integral. It follows directly from the definition of $\Gamma_Q^P$ and the relation \eqref{formula::combinatorial lemma} that, for fixed $Q \in \CF^G(P_0')$,
		\begin{align*}
			\sum_{P \supset Q} \tau_P(X_P^G - T_P^G - T'^G_P) \Gamma_Q^P (X_Q^P  - T_Q^P, T'^P_Q) = \tau_Q(X_Q,T_Q).
		\end{align*}
		Hence we see that $\period^{G',T + T'}(\varphi \otimes \varphi') = \period^{G',T}(\varphi\otimes \varphi')$.
	\end{proof}
	

	\subsection{Periods of truncated automorphic forms}\label{sec::truncated periods}
	
	Another way to approach the regularized integral is through the study of periods of truncated automorphic forms.
	
	As in \cite{Lapid-Rogawski-periods-Galois,Ichino-Yamana-Compositio-Periods,Zydor-periods-AENS}, a key role is played by the Fourier transform of the $\Gamma$ function of Arthur or Zydor. In the last paragraph of \cite[Section 1.4]{Zydor-periods-AENS}, it is asserted there the Fourier transform of the $\Gamma$ function associated to a cone is a polynomial exponential function of the truncation parameter $T$. We do not come up with a proof of this assertion (see Remark \ref{rmk::Zydor's work} below). So we shall make the following assumption on the pair $(G,G')$. For any $P,Q \in \CF^G(P'_0)$ with $P \subset Q$, we have
	\begin{align}\label{formula::condition--truncated-by-Levi}
		A(\overline{\mathfrak{z}_{Q}^+},\overline{\mathfrak{z}_P^+}) \cap \mathfrak{z}_P^Q  \subset A(\overline{\mathfrak{z}_Q^+},\overline{\mathfrak{z}_P^+})^{\vee}.   \tag{$\star$}
	\end{align}
	This is a relative analogue of a standard fact in the geometry of root chambers \cite[Lemma 2.1, (3)]{Zydor-periods-AENS}. 
	
	In view of Lemma \ref{lem::rel parabolics-descent}, if \eqref{formula::condition--truncated-by-Levi} is satisfied for the pair $(G,G')$, then it is also satisfied for $(M,M')$ with $M$ the Levi subgroup of a parabolic subgroup $P \in \CF^G(P'_0)$. Observe that the assumption \eqref{formula::condition--truncated-by-Levi} holds at least in the following two cases:
	\begin{itemize}
		\item[(1)]  $ \mathfrak{a}_0  =  \mathfrak{a}'_{0}$;
		\item[(2)] $(G,G')$ is a symmetric pair. That is, there exists an involution $\theta$ on $G$ such that $G' = G^{\theta}$. 
	\end{itemize}
	For case (1), this is just \cite[Lemma 2.1, (3)]{Zydor-periods-AENS}. For case (2), by \cite[Lemma 2.4]{Helminck-Wang-Involution}, we can choose $A_0$ such that $A_0$ is $\theta$-stable. By \cite[Lemma 3.5]{Helminck-Wang-Involution}, the fixed point group $A_0^{\theta}$ is a maximal $F$-split torus of $G'$ which is taken as our $A_0'$. Thus $\theta$ also acts on $\mathfrak{a}_0$ and $\mathfrak{a}_0^{\theta} = \mathfrak{a}'_0$. The inner product on $\mathfrak{a}_0$ can be taken as $\theta$-invariant. Then \eqref{formula::condition--truncated-by-Levi} follows from \cite[Lemma 2.1, (3)]{Zydor-periods-AENS} and the simple fact that
	\begin{align*}
		A(\overline{\mathfrak{z}_{Q}^+},\overline{\mathfrak{z}_P^+}) \subset A(\overline{\mathfrak{a}_{Q}^+},\overline{\mathfrak{a}_P^+}),\quad \quad \mathfrak{z}_P^Q \subset \mathfrak{a}_P^Q.
	\end{align*}
	We remark that one can find couterexamples for which \eqref{formula::condition--truncated-by-Levi} is not true in some extreme cases where $G'$ is a torus.
	
	\begin{lem}\label{lem::Fourier-Gamma}
		Assume that \eqref{formula::condition--truncated-by-Levi} holds for $(G,G')$. Let $q \in \BC[\mathfrak{z}_P^Q]$ and $\lambda \in \mathfrak{a}'_{0,\BC}$. Let $P \subset Q$ be two parabolic subgroups in $\CF^{G'}(P_0')$ and $T \in \mathfrak{z}^Q$. Then
		\begin{align*}
			\CF(\Gamma_P^Q,T,q,\lambda) := \int_{\mathfrak{z}_P^Q} \Gamma_P^Q(H,T) e^{\langle \lambda,H \rangle} q(H) dH
		\end{align*}
		is a polynomial exponential function of $T \in \mathfrak{z}^Q$, whose exponents may be taken from the set $\{\lambda_R^Q : P \subset R \subset Q\}$. When $\lambda$ is non-degenerate with respect to $A(\overline{\mathfrak{z}_{Q}^+},\overline{\mathfrak{z}_P^+})$, the purely polynomial part of $\CF(\Gamma_P^Q,T,q,\lambda)$ is constant and is given by 
		\begin{align*}
			\int_{\mathfrak{z}_P^Q}^{\sharp} q(H)e^{\langle \lambda,H \rangle} \tau_P^Q(H) dH.
		\end{align*}
	\end{lem}
	\begin{proof}
		The argument is similar to that of \cite[Lemma 2.2]{Arthur-TF-invariant-Annals-81} or \cite[Lemma 4.3]{Zydor-Jacquet-Rallis-Unitary-Canadian-J}. Fix $T$. As $\Gamma_P^Q$ is compactly supported in $H$, $\CF(\Gamma_P^Q,T,q,\lambda)$ is defined by an absolutely convergent integral and is holomorphic at all $\lambda \in \mathfrak{a}'_{0,\BC}$. By the definition of $\Gamma_P^Q$ in \eqref{formula::defn-Gamma}, we have
		\begin{align}\label{formula::Fourier of Gamma of Zydor}
			\CF(\Gamma_P^Q,T,q,\lambda)= \int_{\mathfrak{z}_P^Q} \left(\sum_{\stackrel{R \in \CF^{G'}(P_0')}{P \subset R \subset Q}} \varepsilon_R^Q \tau_P^R(H)\hat{\tau}_R^Q(H-T)   \right) q(H)e^{\langle \lambda,H \rangle} dH.
		\end{align}
	We may take the sum ove $R$ outside the integral if we can find an open set of $\lambda$ such that the integral for each $R$ is absolutely convergent. The existence of such $\lambda$ follow from the claim that 
		\begin{align*}
			(A(\overline{\mathfrak{z}_{R}^+},\overline{\mathfrak{z}_{P}^+}) \cap \mathfrak{z}_P^R)  \times (A(\overline{\mathfrak{z}_{Q}^+},\overline{\mathfrak{z}_{R}^+})^{\vee} \cap \mathfrak{z}_R^Q)  \subset (A(\overline{\mathfrak{z}_{Q}^+},\overline{\mathfrak{z}_{P}^+})^{\vee} \cap \mathfrak{z}_P^Q).
		\end{align*}
		for each $R$ wiht $P \subset R \subset Q$. In fact, by \eqref{formula::condition--truncated-by-Levi}, we have $A(\overline{\mathfrak{z}_{R}^+},\overline{\mathfrak{z}_{P}^+}) \cap \mathfrak{z}_P^R \subset A(\overline{\mathfrak{z}_{R}^+},\overline{\mathfrak{z}_{P}^+})^{\vee} \cap \mathfrak{z}_P^R$, and the latter is contained in $A(\overline{\mathfrak{z}_{Q}^+},\overline{\mathfrak{z}_{P}^+})^{\vee} \cap \mathfrak{z}_P^Q$ by definition. The $1$-dimensional faces of $A(\overline{\mathfrak{z}_{Q}^+},\overline{\mathfrak{z}_{R}^+})^{\vee} \cap \mathfrak{z}_R^Q$ are of the form 
		$A(\overline{\mathfrak{z}_{S}^+},\overline{\mathfrak{z}_{R}^+})^{\vee} \cap \mathfrak{z}_R^S = A(\overline{\mathfrak{z}_{S}^+},\overline{\mathfrak{z}_{R}^+}) \cap \mathfrak{z}_R^S$ with $R \subset S$ and $\dim \mathfrak{z}_R^S = 1$. Note that $(A(\overline{\mathfrak{z}_{S}^+},\overline{\mathfrak{z}_{R}^+}) \cap \mathfrak{z}_R^S)$ is a $1$-dimensional face of $A(\overline{\mathfrak{z}_{S}^+},\overline{\mathfrak{z}_{P}^+}) \cap \mathfrak{z}_P^S$ which is contained in $A(\overline{\mathfrak{z}_{Q}^+},\overline{\mathfrak{z}_{P}^+})^{\vee} \cap \mathfrak{z}_P^Q$ by \eqref{formula::condition--truncated-by-Levi} as above. So when $\lambda$ is negative with respect to $A(\overline{\mathfrak{z}_{Q}^+},\overline{\mathfrak{z}_{P}^+})^{\vee}$ we have 
		\begin{align*}
			\CF(\Gamma_P^Q,T,q,\lambda)=\sum_R \varepsilon_R^Q \int_{\mathfrak{z}_P^Q} \tau_P^R(H^R) \hat{\tau}_R^Q (H_R - T_R) q(H)e^{\langle \lambda,H\rangle} dH.
		\end{align*}
		The first statement follows easily. For a general $\lambda$, we use the similar analysis in \cite{Arthur-TF-invariant-Annals-81} or \cite{Zydor-Jacquet-Rallis-Unitary-Canadian-J}. We omit the details. When $\lambda$ is non-degenerate with respect to $A(\overline{\mathfrak{z}_{Q}^+},\overline{\mathfrak{z}_{P}^+})$, the restriction of $\lambda$ to $\mathfrak{z}_R^Q$ is nonzero unless $R=Q$. Thus the second statement follows.
	\end{proof}
	\begin{rem}\label{rmk::Zydor's work}
		In the general case we do not know how to take sum over $R$ outside the integral in \eqref{formula::Fourier of Gamma of Zydor}. We point out that Zydor's results on regularized periods depend crucially on the properties of the Fourier transform of his $\Gamma$ functions \cite[Theorem 4.1]{Zydor-periods-AENS}.
	\end{rem}

	\begin{prop}
		Assume that \eqref{formula::condition--truncated-by-Levi} holds for $(G,G')$. 
		
		\textup{(i)}	For $\varphi \in \CA(G)$ and $\varphi' \in \CA(G')$, the function
		\begin{align*}
			T \mapsto \int_{[G']^{1,G}} \Lambda^T \varphi(g) \varphi'(g) dg
		\end{align*}
		defined for $T$ sufficiently positive is a polynomial exponential $\sum_{\lambda} p_{\lambda}(T) e^{\langle \lambda,T \rangle}$. The exponents may be taken from the set 
		\begin{align*}
			\mathop{\cup}\limits_{\stackrel{P,Q \in \CF^{G'}(P_0')}{P \subset Q}}  \{ (\lambda + \lambda' + \rho_{P'} + \underline{\rho}_P)_Q^G \ |\ \lambda \in \mathcal{E}_P(G)',\ \lambda' \in \mathcal{E}_{P'}(G') \}.
		\end{align*}
		
		\textup{(ii)} If $(\varphi, \varphi') \in \mathcal{A}(G \times G')^{\ast}$, then $\period^{G'}(\varphi \otimes \varphi') = p_0(T)$. In particualr, the right hand side is constant.
	\end{prop}
	\begin{proof}
		(i) Integrate the relation \eqref{formula::lambda-T-T'} over $[G']^{1,G}$. We get
		\begin{align*}
			\int_{[G']^{1,G}} & \Lambda^{T + T'} \varphi(g) \varphi'(g) dg \\
			= &\sum_P \int_{P'(F) \backslash G'(\BA)^{1,G}} \Lambda^{T,P}\varphi(g) \varphi'(g) \Gamma_P (H_{0'}(g)_P^G - T_P^G, T'^G_P ) dg  \\
			= &\sum_P \int_{P'(F) \backslash G'(\BA)^{1,G}} \Lambda^{T,P}\varphi(g) \varphi'_{P'}(g) \Gamma_P (H_{0'}(g)_P^G - T_P^G, T'^G_P ) dg. 
		\end{align*}
		Expand $\varphi_P$ and $\varphi'_{P'}$ as in \eqref{formula::expansion--before regularization}. The last inner integral is equal to 
		\begin{align*}
			\sum_{ij} \int_{k \in K'} \int_{[M']^{1,P}}\Lambda^{T,P} &\phi_i(mk)\phi'_j(mk)  e^{ - \langle 2 \rho_{P'} , H_{P'}(m) \rangle }  dmdk
		\end{align*}
		times
		\begin{align*}
			 \int_{\mathfrak{z}_P^G} q_i(X)q_j'(X) e^{ \langle \lambda_i  + \lambda'_j + \rho_{P'} + \underline{\rho}_P , X \rangle}  \Gamma_P (X- T_P^G, T'^G_P )  dX. 
		\end{align*} 
		By Lemma \ref{lem::Fourier-Gamma}, the last integral is a polynomial exponential function in $T'$ whose exponents are $\{ (\lambda_i + \lambda'_j + \rho_{P'} + \underline{\rho}_P)_Q^G \}_{Q \supset P}$. This proves \textup{(i)}.
		
		(ii)\ Applying \textup{(i)} to the pair $(M,M')$ we know from definition that each $\period^{G',T}_P(\varphi \otimes \varphi')$ is a polynomial exponential function in $T$. By assumption, the purely polynomial part does not occur unless $P = G$. Since the regularized period is a constant independant of $T$, it is equal to the purely polynomial part of the term $P = G$, which is just $p_0(T)$. 
	\end{proof}
		
	Let $P \in \CF^G(P'_0)$, $\varphi \in \mathcal{A}_P(G)$ and $\varphi' \in \mathcal{A}_{P'}(G')$. We may generalize the construction of Section \ref{sec::regularization of periods} to define the regularized integral
	\begin{align}\label{formula::regularized--Levi}
		\int^*_{P'(F) \backslash G'(\BA)^{1,G}} \varphi (g) \varphi'(g) \tau(H_{0'}(g)^G_P - T^G_P) dg,
	\end{align}
	where $\tau$ is the characteristic function of a cone in $\mathfrak{z}_P^G$, as
	\begin{align*}
		\int_{K'} \int_{\mathfrak{z}_P^G}^{\sharp} \left( \int_{M'(F) \backslash M'(\BA)^{1,M}}^*  \varphi(e^X mk) \varphi'(e^Xmk) e^{- \langle 2\rho_{P'}, H_{P'}(m) \rangle} dm \right) e^{-2\langle \rho_{P'}, X \rangle} \tau(X - T_P^G) dX dk.
	\end{align*}
	For $\tau = \hat{\tau}_P$ this is well-defined provide that the following two conditions are satisfied:
	\begin{equation*}
		\langle \lambda  + \lambda' + \rho_{Q'} + \underline{\rho}_Q, \,\mathfrak{z}_R^P \rangle \neq  0.  \tag{$1^*$}
	\end{equation*}
	for any $Q \subset R \subset P$ in $\CF^G(P'_0)$ with $\dim \mathfrak{z}_R^P = 1$, $\lambda \in \CE_Q(\varphi)$ and $\lambda' \in \CE_{Q'}(\varphi')$.
	\begin{equation*}
		\langle \lambda  + \lambda' + \rho_{P'} + \underline{\rho}_P, \,\mathfrak{z}_P^R \rangle \neq 0 \tag{$2^*$}
	\end{equation*}
	for any $P \subset R$ in $\CF^G(P'_0)$ with $\dim \mathfrak{z}_P^R = 1$, $\lambda \in \CE_P (\varphi)$ and $\lambda' \in \CE_{P'}(\varphi')$.
	
	Let $\CA(G \times G')^{**}$ be the space of pairs $(\varphi,\varphi') \in \CA(G)\oplus \CA(G')$ which satisfy $(1^*)$ (and then also $(2^*)$) for all $P \in \CF^G(P'_0)$.

	\begin{prop}\label{prop::TruncatedPeriod-Main}
		Let $(G,G')$ be such that the assumption \eqref{formula::condition--truncated-by-Levi} holds. If $(\varphi,\varphi') \in \mathcal{A}(G \times G')^{**}$, then
		\begin{align*}
			\int_{[G']^{1,G}} \Lambda^T \varphi ( g ) & \varphi'(g) dg \\
			&= \sum_{P\in \mathcal{F}^G(P_0')} \varepsilon_P^G \int_{P'(F) \backslash G'(\BA)^{1,G}}^* \varphi_P (g) \varphi'_{P'} (g) \hat{\tau}_P (H_{0'}(g)^G - T^G) dg.
		\end{align*}
	\end{prop}
	\begin{proof}
		The argument parallels that in the proof of \cite[Theorem 10]{J-L-R-periods-JAMS}. By induction on the rank of $G$, we may assume that the theorem holds for pairs $(M,M')$ where $M$ is the Levi subgroup of a proper parabolic subgroup $P$ in $\CF^G(P'_0)$. We will show below that the induction hypothesis implies that
		\begin{align}\label{formula::induction hypothesis-I}
			\int_{P'(F) \backslash G'(\BA)^{1,G}}^{\sharp} \Lambda^{T,P}\varphi_P (g) \varphi'_{P'} (g) \tau_P (H_{0'}(g)_P - T_P) dg
		\end{align}
		is equal to 
		\begin{align}\label{formula::induction hypothesis-II}
			\sum_{\stackrel{Q \in \mathcal{F}^G(P_0')}{Q \subset P}} \varepsilon_Q^P \int_{Q'(F) \backslash G'(\BA)^{1,G}}^* \varphi_Q(g)  \varphi'_{Q'} (g)  \hat{\tau}_Q^P (H_{0'}(g)^P - T^P ) \tau_P(H_{0'}(g)_P - T_P) dg.
		\end{align}
		Assuming this, we sum over $P$ to write
		\begin{align}
			\int^*_{G'(F) \backslash G(\BA)^{1,G}} \varphi(g)\varphi'(g) dg - \int_{G'(F) \backslash G(\BA)^{1,G}} \Lambda^T \varphi(g)\varphi'(g) dg
		\end{align}
		as 
		\begin{align}\label{formula::trunc=sum-Levi--4}
			\sum_{\stackrel{Q \subset P}{P \neq G}}  \varepsilon_Q^P \int_{Q'(F) \backslash G'(\BA)^{1,G}}^* \varphi_Q(g)  \varphi'_{Q'} (g)  \hat{\tau}_Q^P (H_{0'}(g)^P - T^P) \tau_P(H_{0'}(g)_P - T_P) dg.
		\end{align}
		As $Q \ne G$, it follows from the relation \eqref{formula::combinatorial lemma} that
		\begin{align*}
			\sum_{\stackrel{Q \subset P}{P \neq G}} \varepsilon_Q^P \hat{\tau}_Q^P (H_{0'}(g)^P - T^P) \tau_P(H_{0'}(g)_P - T_P) =  - \varepsilon_Q^G \hat{\tau}_Q (H_{0'}(g)^G - T^G).
		\end{align*}
		Hence the theorem will follow if we can check that the summation can be taken inside the integral in \eqref{formula::trunc=sum-Levi--4}.
		
		By our assumption \eqref{formula::condition--truncated-by-Levi}, we have 
		\begin{align*}
			(A(\overline{\mathfrak{z}_P^+},\overline{\mathfrak{z}_Q^+})^{\vee} \cap \mathfrak{z}_Q^P)  \times (A(\overline{\mathfrak{z}_G^+},\overline{\mathfrak{z}_P^+}) \cap \mathfrak{z}_P^G) \subset (A(\overline{\mathfrak{z}_G^+},\overline{\mathfrak{z}_Q^+})^{\vee} \cap \mathfrak{z}_Q^G)
 		\end{align*}
		for all $P$ containing $Q$. Let $\lambda \in \CE_Q(\varphi)$, $\lambda' \in \CE_{Q'}(\varphi')$ and $P$ containing $Q$. By our hypothesis,
		\begin{align*}
			\langle \lambda  + \lambda' + \rho_{Q'} + \underline{\rho}_Q, \,\mathfrak{z}_Q^R \rangle \neq  0
		\end{align*}
		for any $R$ with $Q \subset R \subset P$ and $\dim \mathfrak{z}_Q^R = 1$, and
		\begin{align*}
			\langle \lambda  + \lambda' + \rho_{Q'} + \underline{\rho}_Q, \,\mathfrak{z}_S^G \rangle \neq  0
		\end{align*}
		for any $S$ with $P \subset S \subset G$ and $\dim \mathfrak{z}_S^G = 1$. Thus $\lambda + \lambda' + \rho_{ Q'} + \underline{\rho}_Q$ is non-degenerate with respect to the cone 
		\begin{align*}
			(A(\overline{\mathfrak{z}_P^+},\overline{\mathfrak{z}_Q^+})^{\vee} \cap \mathfrak{z}_Q^P)  \times (A(\overline{\mathfrak{z}_G^+},\overline{\mathfrak{z}_P^+}) \cap \mathfrak{z}_P^G).
		\end{align*}
		We can then apply Lemma \ref{lem::Lem 6} and the definition of regularized integrals to justify the exchange of summation and integration as mentioned above.

		It remains to prove the equality of \eqref{formula::induction hypothesis-I} and \eqref{formula::induction hypothesis-II}. Expand $\varphi_P$ and $\varphi'_{P'}$ as in \eqref{formula::expansion--before regularization}:
		\begin{align*}
			\varphi_P(u e^X mk) &= \sum_i q_i(X) e^{\langle \lambda_i + \rho_P,X \rangle}\phi_i(mk)  \\
			\varphi'_{P'} (u e^X mk) &= \sum_j q'_j(X) e^{\langle \lambda'_j + \rho_{P'},X \rangle}  \phi'_j(mk)
		\end{align*}
		with $\phi_i \in \CA_P(G)$ and $\phi_j' \in \CA_{P'}(G')$. Then \eqref{formula::induction hypothesis-I} is equal to the integral over $k \in K'$ and sum over $i,j$ of
		\begin{align*}
			\int_{[M']^{1,M}}   \Lambda^{T,P}\phi_i (mk) &  \phi'_j(mk) e^{-\langle 2\rho_{P'}, H_{P'}(m) \rangle} dm  \\
			&\cdot \int_{\mathfrak{z}_P^G}^{\sharp} q_i(X)q'_j(X) e^{\langle \lambda_i + \lambda'_j + \rho_{P'} + \underline{\rho}_P, X \rangle} \tau_P(X - T_P) dX.
		\end{align*}
		By Lemma \ref{lem::rel parabolics-descent} we can replace $\Lambda^{T,P}$ by $\Lambda^{T,M}$, the mixed truncation with respect to $M$. By induction hypothesis, the first integral over $m \in [M']^{1,G}$ is equal to the sum over $Q \subset P$ of  
		\begin{align}\label{formula::trunc=sum-Levi--5}
			\varepsilon_Q^P \int^{\ast}_{(M'\cap Q')(F) \backslash M'(\BA)^{1,M}}   (\phi_i)_Q (mk)  (\phi'_j)_{Q'}(mk) & e^{-\langle  2\rho_{P'}, H_{P'}(m) \rangle }   \\
			&\hat{\tau}_{M \cap Q}^M ( H_{0'}(m)^M - T^M) dm. \nonumber
		\end{align}		
		We then expand $(\phi_i)_Q$ and $(\phi_j')_{Q'}$ again as in the proof of Proposition \ref{prop::period-well defined}.
		\begin{align*}
			\begin{aligned}
				(\phi_i)_Q ( v e^Y lk'k) &= \sum_{s} q_{is} (Y)e^{\langle \lambda_{is} + \rho_Q,Y \rangle}    \phi_{is} (lk'k) \\
				(\phi'_i)_{Q'} ( v e^Y l k'k)  &= \sum_{t} q'_{jt} (Y)e^{\langle \lambda'_{jt} + \rho_{Q'}, Y \rangle}   \phi'_{jt} (lk'k)
			\end{aligned}
		\end{align*}
		with $m = v e^Y l k'$ an Iwasawa decomposition of $m$ with respect to $M \cap Q$. Then \eqref{formula::induction hypothesis-I} is equal to the integral over $k \in K'$, sum over $Q \subset P$ and over $i,j,s,t$ of the product of
		\begin{align*}
			\varepsilon_Q^P \int_{[L']^{1,Q}}^{\ast}  \phi_{is}(lk) \phi'_{jt}(lk) e^{-\langle  2\rho_{Q'},  H_{Q'}(l) \rangle } dl
	    \end{align*} 
        with
        \begin{align*}
			\int^{\sharp}_{\mathfrak{z}_Q^P}  q_{is}(Y)q'_{jt}(Y)   & e^{\langle \lambda_{is} + \lambda'_{jt} + \rho_{Q'} + 2 \underline{\rho}_Q, Y \rangle} \hat{\tau}_Q^P (Y - T^P) dY  \\
			 &\cdot \int^{\sharp}_{\mathfrak{z}_P^G} q_i(X)q_j'(X) e^{\langle \lambda_i + \lambda'_j + \rho_{P'} + \underline{\rho}_P, X \rangle} \tau_P(X - T_P) dX.
		\end{align*}
	   Here we have used the relation \eqref{formula::rho-relation} and the relation
     	\begin{align*}
		    \hat{\tau}_{M \cap Q}^M ( H_{0'}(m)^M - T^M) = \hat{\tau}_{Q}^P( H_{0'}(m)^P - T^P). 
    	\end{align*}
		By \eqref{formula::combine-integrals-I}, we can combine the two $\sharp$-integrals into one $\sharp$-integral over $\mathfrak{z}_Q^G$. Then the equality of \eqref{formula::induction hypothesis-I} and \eqref{formula::induction hypothesis-II} follows directly from the definition of the regularized integral over $Q'(F) \backslash G'(\BA)^{1,G}$.
	\end{proof}

	\section{Global linear periods: The cuspidal case}\label{sec::linear-periods-cuspidal}
	
	\subsection{Notations}\label{sec::cusp linear periods-notation}
	We fix some notations that will be used from now on. For any positive integer $n$, denote by $G_n$ the general linear group $\GL_n$ defined over $F$, and by $Z_n$ its center. Denote by $B_n$, $T_n$ and $N_n$ the subgroup of upper triangular matrices, the subgroup of diagonal matrices and the subgroup of upper triangular unipotent matrices in $G_n$ respectively.
	
	By a partition of a positive integer $n$, we mean a tuple $\bar{n}= (n_1,\cdots,n_t)$ of positive integers such that $n = n_1 + \cdots + n_t$. To such a partition, we associate a parabolic subgroup $P_{\bar{n}}$, consisting of elements of $G_n$ of the form
	\begin{align*}
		\begin{pmatrix}
			g_1 & \ast  & \ast \\
			    & \ddots& \ast \\
			    &       & g_t
		\end{pmatrix}
	\end{align*} 
	with $g_i \in G_{n_i}$ for $1 \leqslant i \leqslant t$. We denote by $N_{\bar{n}}$ its unipotent radical and $M_{\bar{n}}$ its standard (containing $T_n$) Levi subgroup consisting of matrices $\textup{diag}(g_1,\cdots,g_t)$ with $g_i \in G_{n_i}$. We denote by $\CP_n = \mathcal{M}_n\mathcal{U}_n$ the parabolic subgroup of $G_n$ associated to the partition $(n-1,1)$, and by $\mathcal{U}_n$ its unipotent radical. We denote by $P_n$ the mirabolic subgroup of $G_n$ consisting of matrices whose last row is $e_n = (0,\cdots,0,1) \in F^n$. Denote
	\begin{align*}
		P_n^{(i)} = \left\{ \begin{pmatrix}
			g & x \\
			   & u
		\end{pmatrix} \ |\ g \in G_{n-i}, u \in N_i, x \in M_{n-i,i}   \right\}
	\end{align*}
	for $1 \leqslant i \leqslant n$. Denote by $U_n^{(i)}$ the unipotent radical of $P_n^{(i)}$. Sometimes we view $G_{n-1}$ as a subgroup of $G_n$ via the map $g \mapsto \textup{diag}(g,1)$.

	We denote by $\wo$ the element of the symmetric group $\mathfrak{S}_{2n}$ naturally embedded in $G_{2n}$ defined by 
	\begin{align*}
		\begin{pmatrix}
			1 & 2 & \cdots & n-1 & n & n+1 & n+2 & \cdots & 2n-1 & 2n  \\
			1 & 3 & \cdots & 2n-3 & 2n-1 & 2 & 4 & \cdots & 2n-2 & 2n
		\end{pmatrix}.
	\end{align*}
	Let $G_{2n}' = w_{\star} M_{(n,n)} w_{\star}^{-1}$. If we denote by $\epsilon_{2n}$ the diagonal matrix 
	\begin{align*}
		\left(\begin{smallmatrix}
			1 &  & & & &  \\
			  & -1 & & & &  \\
			  & & 1 & & & \\
			  & & & -1 & & \\
			  & & & & \ddots & \\
			  & & & & & -1
		\end{smallmatrix} \right)
	\end{align*}
	of $G_{2n}$, $G'_{2n}$ is the subgroup of elements fixed by the involution $g \mapsto \epsilon_{2n} g \epsilon_{2n}$. Our main concern in this paper is within the symmetric pair $(G_{2n},G^{\prime}_{2n})$. If $S$ is a subset of $G_{2n}$, then, unless otherwise stated, $S'$ will stand for the intersection of $S$ with $G'$. For $g_1,g_2 \in G_n$, we will write 
	\begin{align}\label{formula::defn-element-G'}
		\iota(g_1,g_2) = w_{\star} \textup{diag}(g_1,g_2)w_{\star}^{-1} \in G'_{2n}. 
	\end{align}
	We then can identify $G_n \times G_{n-1}$ with a subgroup of $G'_{2n}$ via the map $\iota$.  
	
	For $g \in G_n(k)$ where $k$ is a local field or the adele ring of a number field, we sometimes write $|g|$ or $\nu (g)$ for $| \det g|$, the norm of the determinant of $g$.
	
	For $s \in \BC$, we denote by $\mu_s$ the character on $G'_{2n}(\BA)$ defined by
	\begin{align*}
		\mu_s (\iota(g_1,g_2)) = |\det g_1|^s|\det g_2|^{-s}.
	\end{align*}
	
	We fix a nontrivial additive character $\psi$ of $F \backslash \BA$. Denote by $\CS(\BA^n)$ the space of Schwartz functions on $\BA^n$. For $\Phi \in \CS(\BA^n)$, write for $\hat{\Phi}$ the Fourier transform of $\Phi$ with respect to a $\psi$-self-dual Haar measure, given by
	\begin{align*}
		\hat{\Phi}(X)  = \int_{\BA^n} \Phi(Y) \psi( - \Tr(X {}^tY)) dY.
	\end{align*}

	\subsection{Global linear periods}
	
	Let $\pi = \otimes_v \pi_v$ be a cuspidal automorphic representation of $G_{2n}(\BA)$ with trivial central character. Then we have the global (twisted) linear period
	\begin{align}\label{formula::defn-period-integral}
		\period(\varphi,s_0) := \int_{G_{2n}'(F) \backslash G'_{2n}(\BA)^{1,G}} \varphi(g) \mu_{s_0}(g) dg, \quad \varphi \in \pi,\ s_0 \in \BC.
	\end{align}
	The integral converges absolutely by \cite[Proposition 1]{Ash-Ginzburg-Rallis-vanishing-pair}. In the literature, the integral in the definition of the linear period is often taken over the quotient 
	\begin{align*}
		Z_{2n}(\BA)G'_{2n}(F) \backslash G'_{2n}(\BA).	
	\end{align*}
	Observe that $G'_{2n}(\BA)^{1,G} = G'_{2n}(\BA) \cap G_{2n}(\BA)^1$. So we have a surjective map
	\begin{align*}
		G_{2n}'(F) \backslash G'_{2n}(\BA)^{1,G}  \ra Z_{2n}(\BA)G'_{2n}(F) \backslash G'_{2n}(\BA)
	\end{align*}
	with each fiber isomorphic to $F^{\times} \backslash \BA^{\times}$. The linear period is closely related to the Bump-Friedberg $L$-function introduced in \cite{Bump-Friedberg-Exterior-L} (see \cite{Matringe-specialisation-BF-L}).
	
 	\subsection{Rankin-Selberg theory of Bump-Friedberg $L$-functions} 
 	
 	We first recall the mirabolic Eisenstein series used in the study of Bump-Friedberg $L$-functions. For $\Phi \in \CS(\BA^n)$ and $s \in \BC$, the integral
	\begin{align*}
		F(g,\Phi;s) = |\det g|^s \int_{\BA^{\times}} \Phi(e_n a g)|a|^{ns} d^{\times} a,\quad g \in G_n(\BA),
	\end{align*}
	converges for $\Re s > 1/n$. Consider the Eisenstein series $\tilde{E}(\Phi,s)$ on $G_n(\BA)$ defined by
	\begin{align}\label{formula::defn--mirabolic Eisenstein}
		\tilde{E}(g,\Phi;s) = \sum_{\gamma \in \CP_n(F) \backslash G_n(F)} F(\gamma g, \Phi;s).
	\end{align}
	We recall the following fact from \cite[Lemma 4.2]{J-S-Euler-product-I} and \cite[Proposition 3.1]{ZhangWei-automorphic-period-JAMS-2014}.
	\begin{lem}\label{lem::epstein-Eisenstein-pole}
		The series $\tilde{E}(\Phi,s)$ converges absolutely for $\Re s >1$ and admits a meromorphic extension to $\BC$. It has a simple pole at $s = 1$ with residue $\frac{\hat{\Phi}(0)}{n}v_F$, where $v_F = \textup{vol}(F^{\times} \backslash \BA^1)$.
	\end{lem}
	
	Following \cite{Matringe-specialisation-BF-L}, for $\Phi \in \CS(\BA^n)$ and $s_1,s \in \BC$, we set
	\begin{align*}
		f(\iota(g_1,g_2),\Phi;s_1,s) = | g_1 |^{s_1 + s -1/2} |  g_2 |^{s_1 - s + 1/2} \int_{\BA^{\times}} \Phi ( e_n a g_2) |a|^{2ns_1} d^{\times}a 
	\end{align*}
	for $\iota(g_1,g_2) \in G_{2n}'(\BA)$. The integral converges for $\Re(s_1) > 1/2n$. Observe that $f(\Phi ; s_1,s)$ is left $Z_{2n}(\BA)\CP'_{2n}(F)$-invariant. Consider the Eisenstein series $E(\Phi;s_1,s)$ on $G'_{2n}(\BA)$ defined by
	\begin{align}\label{formula::Eisenstein series--linearperiod}
		E(g,\Phi; s_1,s) = \sum_{\gamma \in \CP'_{2n}(F) \backslash G'_{2n}(F)} f(\gamma g, \Phi;s_1,s).
	\end{align}	 
	Note that we have
	\begin{align}\label{formula::connection-Eisenstein-Mirabolic}
		E(\iota(g_1,g_2),\Phi;s_1,s) = |g_1|^{s_1 + s - 1/2} | g_2 |^{-s_1 - s + 1/2} \tilde{E}(g_2,\Phi;2s_1).
	\end{align}
	So it follows from Lemma \ref{lem::epstein-Eisenstein-pole} that
	
	\begin{lem}\label{lem::pole--mirabolic Eisenstein}
		The series $E(\Phi;s_1,s)$ converges absolutely for $\Re s_1 > \frac{1}{2}$ and admits a meromorphic extension to $\BC \times \BC$. It has simple poles along $s_1 = \frac{1}{2}$ with residues 
		\begin{align*}
			\res_{s_1 = \frac{1}{2}}E(\iota(g_1,g_2),\Phi;s_1,s)  =  v_F \frac{\hat{\Phi}(0)}{2n} |\det g_1|^s |\det g_2|^{-s}.
		\end{align*}
	\end{lem}

	For $\varphi \in \pi$, let
	\begin{align*}
		W_{\varphi}^{\psi} (g) = \int_{N_{2n}(F) \backslash N_{2n}(\BA)} \varphi (ng) \overline{\psi(n)} dn
	\end{align*}
	be the Whittaker function associated to $\varphi$ (and $\psi$). The following result summarizes what Matringe proved in \cite{Matringe-specialisation-BF-L}.
	\begin{prop}\label{prop::Bump-F}
		Let $\pi$ be an irreducible cuspidal automorphic representation of $G_{2n}(\BA)$ with trivial central character. Assume $\Re s_0 \geqslant 0$. Then
		\begin{align}\label{formula::Bump-F}
			&	\int_{[G'_{2n}]^{1,G_{2n}}} \varphi(g) E(g,\Phi;s_1,s_0)dg \\ 
			&= v_F \int_{(N_n(\BA) \backslash G_n(\BA))^2 } W_{\varphi}^{\psi}(\iota(g_1,g_2)) \Phi(e_n g_2) | g_1|^{s_1 + s_0 - 1/2} | g_2|^{ s_1 -(s_0 - 1/2)} dg_1dg_2.  \nonumber
		\end{align}
	for $\varphi \in \pi$ and $\Phi \in \CS(\BA)^n$, where the right hand side converges absolutely for $\Re s_1$ sufficiently large. If $W_{\varphi}^{\psi} (g) = \prod_v W_v(g_v)$ and $\Phi = \otimes_v \Phi_v$, the right hand side decomposes into a product of local integrals
	\begin{align*}
		\prod_v \Psi(s_1,W_{\varphi,v},\Phi_v,\mu_{s_0}),
	\end{align*}
	where the local integral is defined similarly as the right hand side of \eqref{formula::Bump-F}. The infinite product of the local integrals converges absolutely when $\Re s_1$ sufficiently large. Each local integral converges for $\Re s_1 \geqslant \frac{1}{2}-\varepsilon$ with $\varepsilon$ a positve real number. For unramified data of $(W_{\varphi,v},\Phi_v)$, one has
	\begin{align}\label{formula::unramified-computation}
		\Psi(s_1,W_{\varphi,v},\Phi_v,\mu_{s_0}) = L(s_1 + s_0,\pi_v) L (2s_1,\pi_v,\Lambda^2).
	\end{align}
	Here $L(s,\pi_v)$ and $L(s,\pi_v,\Lambda^2)$ are the corresponding $L$-functions of the Langlands parameter of $\pi_v$ (\cite{Matringe-BF+L-function}).
	\end{prop}
	
	Let $S$ be a finite set of places of $F$ large enough such that for $v$ outside $S$, $\pi_v$ is unramified. As in \cite{Matringe-specialisation-BF-L}, we define an variant of the partial Bump-Friedberg $L$-function as
	\begin{align*}
		L^{\textup{lin},S} (s,\pi,\mu_{s_0}) = \prod_{v \notin S} L (s + s_0,\pi_v) L(2s,\pi_v,\Lambda^2).
	\end{align*}
	It follows from the unramified computation \eqref{formula::unramified-computation} and the unfolding identity \eqref{formula::Bump-F} that $L^{\textup{lin},S}(s,\pi,\mu_{s_0})$ is a meromorphic function of $s$.
	
	Following the argument in \cite[p.184]{Gelbart-Jacquet-Rogawski-Generic+U3}, we will give an explicit factorization of the linear periods (cf. \cite[Proposition 3.2]{ZhangWei-automorphic-period-JAMS-2014} for the case of Flicker-Rallis periods). For a place $v$ of $F$, define the local linear form $\vartheta_v$ as follows: For $W_v \in \mathcal{W}(\pi_v,\psi_v)$,
	\begin{align*}
		\vartheta_v(W_v) = \int_{N_n(F_v)\backslash G_v(F_v)} \int_{N_{n-1}(F_v) \backslash G_{n-1}(F_v)} W_v(\iota(g_1,\left(\begin{smallmatrix}
			g_2 &  \\
			    & 1
		\end{smallmatrix}\right))|g_1|^{s_0} |g_2|^{-s_0} dg_1dg_2.
	\end{align*}
	For $v$ such that $\pi_v$ is unramified, we define a normalized linear form
	\begin{align*}
		\vartheta^{\natural}_v(W_v) = \frac{\vartheta_v(W_v)}{L(\frac{1}{2} + s_0,\pi_v)L(1,\pi_v,\Lambda^2)}.
	\end{align*}
	Note that the local factor $L(s+s_0,\pi_v)L(2s,\pi_v,\Lambda^2)$ has no pole or zero at $s = 1/2$ for a unitary unramified generic $\pi_v$.
	
	\begin{prop}\label{prop::explicit decomposition}
		Assume $\Re s_0 \geqslant 0$. Let $S$ be a finite set of places of $F$ sufficiently large such that outside $S$, $\pi_v$ is unramified. We have an explicit decomposition
		\begin{align}\label{formula::explicit decomposition}
			\period(\varphi,s_0) =  2n \res_{s= \frac{1}{2}} L^{\textup{lin},S}(s,\pi,\mu_{s_0}) \prod_{v \in S} \vartheta (W_v) \prod_{v \notin S} \vartheta^{\natural}(W_v),
		\end{align}
		where $W_{\varphi} = \prod_v W_v$. The right hand side does not depend on the choice of $S$.
	\end{prop}
	\begin{proof}
		We assume that the residue of $L^{\textup{lin},S}(s,\pi,\mu_{s_0})$ at $s = 1/2$ is nonzero. If not so, $\period(\varphi,s_0)  = 0$ for all $\varphi \in \pi$ by \cite[Theorem 4.5]{Matringe-specialisation-BF-L} and \eqref{formula::explicit decomposition} holds trivially. We then show that for each place $v$ of $F$, we have
		\begin{align}\label{formula::explicit factorization--I}
			\Psi(\tfrac{1}{2},W_{v},\Phi_v,\mu_{s_0})  =  \vartheta_v(W_{v}) \hat{\Phi}_v(0).
		\end{align}
		For this, we choose a finite set $S'$ of places of $F$ sufficiently large such that $S'$ contains the given place $v$ and all archimedean places, and that outside $S$, $\pi_w$ is unramified, $W_w$ is $K_w$-invariant and equals to $1$ on $K_w$ ($\psi_w$ has conductor $\mathcal{O}_w$), and $\Phi_w$ is the characteristic function of $\mathcal{O}_w^n$. Let $\varphi$ and $\Phi$ be pure tensor with prescribed local component at $v$. Taking residue at $s_1 = 1/2$ in \eqref{formula::Bump-F}, we have
		\begin{align*}
			\period(\varphi,s_0) \hat{\Phi}(0) = 2n \cdot \res_{s= \frac{1}{2}} L^{\textup{lin},S'}(\tfrac{1}{2},\pi,\mu_{s_0}) \prod_{w \in S'} \Psi(\tfrac{1}{2},W_w,\Phi_w,\mu_{s_0}).
		\end{align*}
		Here we used the fact that the local integral converges at $s= 1/2$. For the given place $v$, by the integration formula in Lemma \ref{lem::integration formula-cuspidal}, we have
		\begin{align}\label{formula::explicit factorization--II}
			\Psi(\tfrac{1}{2},W_{v},\Phi_v,\mu_{s_0})  = \int_{\Xi_{n,v}} \Phi_v(X) \vartheta_v(\pi_v(\iota(\mathbf{1}_n,n(X)))W_v) |X_n|^{-s} dX.
		\end{align}
		Here $X = (X_1,\cdots,X_n) \in \Xi_{n,v} = F_v^{n-1} \oplus F_v^{\times} \subset F_v^n$, $n(X)$ is the matrix with last row $X$ and $n(X)_{ij} = \delta_{ij}$ when $i < n$, and $dX$ is the restriction to $\Xi_{n,v}$ of the $\psi_v$ self-dual Haar measure on $F_v^n$. For each place $w \in S'$ other than $v$, by \cite[Corollay 3.3, Proposition 3.5]{Matringe-specialisation-BF-L}, choose $W_w$ and $\Phi_w$ appropriately such that $\Psi(1/2,W_w,\Phi_w,\mu_{s_0}) \ne 0$. Varying $\Phi_v$, we get that $ \vartheta_v(\pi_v(\iota(\mathbf{1}_n,n(X)))W_v) |X_n|^{-s}$ is a constant function that takes value $\vartheta_v(W_v)$. Thus \eqref{formula::explicit factorization--I} follows from \eqref{formula::explicit factorization--II}. Our decomposition \eqref{formula::explicit decomposition} follows then from the unramified computation \eqref{formula::unramified-computation}.
	\end{proof}

	\section{Local linear periods}\label{sec::local periods}

	In this section let $F$ be a local field with characteristic $0$. Let $\sigma$ be an irreducible unitary generic representation of $G_r(F)$, and $P$ the standard parabolic subgroup  of $G_{dr}$ of type $(r,\cdots,r)$. The representation 
	$$\Ind_P^G \sigma^{\boxtimes d} \otimes \delta_P^{1/(2r)}$$
	has a unique irreducible quotient, which will be denoted by $\speh(\sigma,d)$. Here $\Ind_P^G$ stands for the normalized parabolic induction. There representations are the local components of irreducible automorphic representation  appeared in the discrete spectrum of $G_{dr}(\BA)$ (see Section \ref{sec::constant term Eisenstein series}).
		
	Assume now that $F$ is a nonarchimedean local field. We recall some notations from the theory of Bernstein-Zelevinsky derivatives in the representation theory of $G_n(F)$. the functors $\Psi^+$, $\Psi^-$, $\Phi^+$ and $\Phi^-$ have been defined in \cite[Section 3.2]{B-Z-I}. Let $\sigma$ be a smooth representation of $P_n$, then $\Phi^+ (\sigma) = \cind_{P_n\CU_{n+1}}^{P_{n+1}} \sigma[1/2] \boxtimes \psi|_{\CU_{n+1}}$ is a smooth representation of $P_{n+1}$, where $\sigma[1/2]$ is the representation of $P_n$ acting on the space of $\sigma$ by $\sigma[1/2](g) = |\det g|^{1/2} \sigma(g)$. Also, here and in the after, we use $\cind$ and $\ind$ for the (unnormalized) compactly induction and induction respectively in the sense of \cite[Section 2.22, Section 2.21]{B-Z-Survery}. $\Phi^-(\sigma)$ is defined to be the normalized Jacquet module of $\sigma$ with repect to $\CU_n$ and the character $\psi|_{\CU_n}$, regarded as a representation of $P_{n-1}$; $\Psi^-(\sigma)$ is the normalized Jacquet module of $\sigma$ with respect to $\CU_n$ and the trivial character, regarded as a representation of $G_{n-1}$. Let $\tau$ be a smooth representation of $G_n$, then $\Psi(\tau) $ is the representation of $P_{n+1}$ on the space of $\tau$ such that $\CU_{n+1}$ acts trivially and $G_n$ acts by $\tau[1/2]$.
	
	For $1 \leqslant i \leqslant n$ and a smooth representation $\pi$ of $G_n(F)$ of finite length, define $\pi^{(i)} = \Psi^-(\Phi^-)^{i-1}(\pi|_{P_n})$ to be the $i$-th derivative of $\pi$. If there exists $ 1 \leqslant h \leqslant n$ such that $\pi^{(h)} \ne 0$ but $\pi^{(i)} = 0$ for all $i > h$, then $\pi^{(h)}$ is called the highest derivative of $\pi$.
	
	\begin{lem}
		Let $\sigma$ be an irreducible unitary generic representation of $G_r(F)$. The highest derivative of $\speh(\sigma,d)$ is isomorphic to $\speh(\sigma,d-1)[-1/2]$.
	\end{lem}
	For a proof, see \cite[Lemma 5.3]{Yamana-Residue-Periods-JFA}. We recall the following result from \cite[Proposition 3.7]{Yangchang-linear-periods}.
	\begin{lem}\label{lem::mirabolic machine-local}
		Let $\sigma$ be a representation of $P_{n-1}$ and $a,b\in \BC$. We have
		\begin{align*}
			\Hom_{M_{p,q}(F) \cap P_n(F)} (\Phi^+\sigma,\nu^a \boxtimes \nu^b) \cong \Hom_{M_{q-1,p}(F) \cap P_{n-1}(F)} (\sigma,\nu^{b - 1/2} \boxtimes \nu^{a + 1/2}),
		\end{align*}
		where $M_{p,q}$ is the standard Levi subgroup of $G_n$ of type $(p,q)$ and $\nu$ is the character defined by $\nu(g) = |\det g|$.
	\end{lem}
		
	Denote by $\mathcal{S}(F^q)$ the space of Schwartz-Bruhat functions on $F^q$ and by $\mathcal{S}_0(F^q)$ the subspace of $S(F^q)$ consisting of those functions that vanish at $0$. We define an action of $M_{p,q}(F)$ on $\CS(F^q)$ by $R(\textup{diag}(g_1,g_2)) \Phi(x) = \Phi (x g_2)$.
	
	\begin{thm}\label{thm::generic uniqueness}
		Let $\pi$ be a smooth representation of $G_n(F)$ of finite length. Let $p$ and $q$ be two positive integers with $p + q =n$. Then, for all pairs $(a,b) \in \BC^2$ outside a finite union of hyperplanes,  
		\begin{align*}
			\dim_{\BC} \Hom_{M_{p,q}(F)} (\pi \otimes \mathcal{S}(F^q),\nu^a\boxtimes \nu^b) \leqslant \dim_{\BC} \pi^{(n)}.
		\end{align*}
	\end{thm}
	\begin{proof}
		The arguments are similar to those in Theorem 1 of \cite{Kable-Asai-L-AJM}. The space 
		\begin{align*}
			\Hom_{M_{p,q}(F)}(\pi,\nu^a \boxtimes \nu^b)
		\end{align*}
		is zero unless one of the central characters of the irreducible subquotients of $\pi$ equals to $\nu^{pa+qb}$. Hence for $(a,b)$ outside a finite union of hyperplanes,  
		\begin{align*}
			\Hom_{M_{p,q}(F)}(\pi \otimes \CS(F^q),\nu^a \boxtimes \nu^b)
		\end{align*}
		is a subspace of 
		\begin{align*}
			\Hom_{M_{p,q}(F)}(\pi \otimes \CS_0(F^q), \nu^a \boxtimes \nu^b).
		\end{align*}
		It is enough to bound the dimension of the latter space above as in the theorem. We have
		\begin{align*}
			\Hom_{M_{p,q}(F)}(\pi \otimes \CS_0(F^q), \nu^a \boxtimes \nu^b) & \cong \Hom_{M_{p,q}(F)} (\pi \otimes \cind_{M_{p,q}(F)\cap P_n(F)}^{M_{p,q}(F)} \mathbf{1}, \nu^a \boxtimes \nu^b)  \\
			& \cong \Hom_{M_{p,q}(F)} (\pi, \ind_{M_{p,q}(F)\cap P_n(F)}^{M_{p,q}(F)}   \nu^a \boxtimes \nu^{b+1} )    \\
			& \cong \Hom_{M_{p,q}(F) \cap P_n(F)} (\pi,\nu^a \boxtimes \nu^{b+1}).
		\end{align*}
		The last isomorphism is due to the Frobenius duality theorem \cite[Theorem 2.28]{B-Z-Survery}. According to \cite[Section 3.5]{B-Z-I}, the restriction of $\pi$ to $P_n(F)$ has a filtration that has composition factors $(\Phi^+)^{i-1} \Psi^+ (\pi^{(i)})$, $i = 1,\cdots,n$. We apply Lemma \ref{lem::mirabolic machine-local} repeatedly to analyze the space
		\begin{align}
			\Hom_{M_{p,q}(F) \cap P_n(F)} ((\Phi^+)^{i-1}\Psi^+(\pi^{(i)}),\nu^a \boxtimes \nu^{b+1}).
		\end{align}
		Such an analysis has been performed in \cite[Proposition 6.13]{Yangchang-linear-periods}. We treat only the case of even $i$. The case of odd $i$ is similar. Suppose $i = 2k$. If $p > k-1$ and $q > k$, we have
		\begin{align*}
			&\Hom_{M_{p,q}(F) \cap P_n(F)}((\Phi^+)^{i-1} \Psi^+ (\pi^{(i)}),\nu^a \boxtimes \nu^{b+1}) \\
			\cong &\Hom_{M_{q-k,p-k+1} \cap P_{n-i+1}(F)} (\Psi^+(\pi^{(i)}), \nu^{b+1/2} \boxtimes \nu^{a+1/2})  \\
			\cong &\Hom_{M_{q-k,p-k}(F)} (\pi^{(i)},\nu^b \boxtimes \nu^a).
		\end{align*}	
		Hence, by considering central characters, 
		\begin{align*}
			\Hom_{M_{p,q}(F) \cap P_n(F)} ((\Phi^+)^{i-1}\Psi^+(\pi^{(i)}),\nu^a \boxtimes \nu^{b+1})  = 0
		\end{align*}
		for all $(a,b)$ outside a finite union of hyperplanes. If $ p < k$ or $q \leqslant k$, then there exists $i_0 \geqslant 0$ such that	
		\begin{align*}
			&\Hom_{M_{p,q}(F) \cap P_n(F)}((\Phi^+)^{i-1} \Psi^+ (\pi^{(i)}),\nu^a \boxtimes \nu^{b+1})  \\
			\cong &\Hom_{P_{n-i+i_0+1}(F)} ((\Phi^+)^{i_0}\Psi^+(\pi^{(i)}),\nu^{a'}),
		\end{align*}
		with $a'$ depending on $i$, $a$ and $b$. If $i_0 > 0$, then, by the definition of $\Phi^+$ and a variant of the Frobenius duality theorem \cite[Proposition 2.29]{B-Z-Survery}, 
		\begin{align*}
			\Hom_{P_{n-i+i_0+1}(F)} ((\Phi^+)^{i_0}\Psi^+(\pi^{(i)}),\nu^{a'} ) = 0.
		\end{align*}
		If $i_0 = 0$, then
		\begin{align*}
			\Hom_{P_{n-i+1}(F)} (\Psi^+(\pi^{(i)}),\nu^{a'} ) \cong \Hom_{G_{n-i}(F)}(\pi^{(i)},\nu^{a'-1/2} ).
		\end{align*}
		Hence, unless $i = n$, this latter space is zero for all $(a,b)$ outside a finite union of hyperplanes. When $ i_0 = 0$ and $ i =n$, as $G_0$ is the trivial group, we have
		\begin{align*}
			\Hom_{M_{p,q}(F) \cap P_n(F)}((\Phi^+)^{n-1} \Psi^+ (\pi^{(n)}),\nu^a \boxtimes \nu^{b+1})  \cong \Hom_{\BC}(\pi^{(n)},\BC).
		\end{align*}
		This gives the bound as that in the theorem.
	\end{proof}
	
	The analysis in the proof of the preceding theorem yields the following vanishing result.
	\begin{prop}\label{prop::vanishing--unequal case}
		Let $\pi$ be a smooth representation of $G_n(F)$ of finite length. Let $p,q$ be two positive integers with $p + q = n$. Assume $p > q$. Then, for all but finitely many $s \in \BC$,
		\begin{align*}
			\Hom_{M_{p,q}(F)}(\pi, \nu^{qs} \boxtimes \nu^{-ps} )  = 0.
		\end{align*}
	\end{prop}

	We collect here some integration formulas that will be needed elsewhere in this work. The proof of these formulas are all based on the fact that ${}^t\CP_n(F) \mathcal{U}_n(F)$ is a dense open subset of $G_n(F)$. We omit the details here. Here $F$ is allowed to be archimeadean. For any positive integer $k$, set $\Xi^k = F^{k - 1}\oplus F^{\times} \subset F^k$. If $X \in F^{k-1}$ and $z \in F^{\times}$, we denote by $dXdz$ the restriction to its open subset $\Xi^k$ of the Haar measure on $F^k$.
	\begin{lem}\label{lem::integration formula-cuspidal}
		Let $f: G_n(F) \times G_n(F) \ra \BC$ be a continuous function that is left $N_n(F) \times N_n(F)$-invariant. Then
		\begin{align*}
			&\int\limits_{N_n(F) \times N_n(F) \backslash G_n(F) \times G_n(F)}  f(g_1,g_2) | \det g_2 | d(g_1,g_2)  \\
			&= \int_{\Xi^n} \int\limits_{N_n(F) \times N_{n-1}(F) \backslash \GL_n(F) \times \GL_{n-1}(F)} f  \left( g_1, \begin{pmatrix}
				g_2  & 0  \\
				X  & z
			\end{pmatrix}\right)  d(g_1,g_2) dXdz.
		\end{align*}
	\end{lem}
	\begin{lem}\label{lem::integration formula}
		Let $f: G_n(F) \times G_m(F) \ra \BC$ be a continuous function that satisfies
		\begin{align*}
			f(p_1g_1,p_2g_2) = | \det p_1|^r | \det p_2|^{s-1} f(g_1,g_2),
		\end{align*}
		for $p_1 \in P^{(r)}_n(F)$, $g_1 \in G_n(F)$ and $p_2 \in P^{(s)}_m(F)$, $g_2 \in G_m(F)$.
		Then
		\begin{align*}
			&\int\limits_{P^{(r)}_n(F) \times P^{(s)}_m(F) \backslash G_n(F) \times G_m(F)}  f(g_1,g_2) | \det g_2 | d(g_1,g_2)  \\
			&= \int_{\Xi^m} \int\limits_{P^{(r)}_n(F) \times P^{(s-1)}_{m-1}(F) \backslash \GL_n(F) \times \GL_{m-1}(F)} f  \left( g_1, \begin{pmatrix}
				h  & 0  \\
				X  & z
			\end{pmatrix}\right)  d(g_1,h) dXdz.
		\end{align*}
		Similarly, let $f: G_n(F) \times G_m(F) \ra \BC$ be a continuous function that satisfies
		\begin{align*}
			f(p_1g_1,p_2g_2) = | \det p_1 |^{r-1} |\det p_2|^s f(g_1,g_2)
		\end{align*}
		for $p_1 \in P^{(r)}_n(F)$, $g_1 \in G_n(F)$ and $p_2 \in P^{(s)}_m(F)$, $g_2 \in G_m(F)$. Then
		\begin{align*}
			&\int\limits_{P^{(r)}_n(F) \times P^{(s)}_m(F) \backslash G_n(F) \times G_m(F)}  f(g_1,g_2) | \det g_1 | d(g_1,g_2)  \\
			&= \int_{\Xi^n} \int\limits_{P^{(r-1)}_n(F) \times P^{(s)}_{m-1}(F) \backslash G_{n-1}(F) \times G_{m}(F)} f  \left( \begin{pmatrix}
				h  & 0  \\
				X  & z
			\end{pmatrix}, g_2 \right)  d(h,g_2) dXdz.
		\end{align*}
	  Here $dXdz$ is the Haar measure on $F^m$ or $F^n$ restricted to $\Xi^m$ or $\Xi^n$.
	\end{lem}

	\section{Eisenstein series, intertwining operators and multiple-residues}\label{sec::Eisenstein series-const terms}
	
	In this section we will compute the constant terms of square-integrable automorphic forms of $G_{2n}(\BA)$ and the Eisenstein series $E(\Phi,s_1;s)$ given in \eqref{formula::Eisenstein series--linearperiod}.
	
	\subsection{General notations}
	
	Let us retain the notation of Section \ref{sec::notation}. We work with a reductive group $G$ over $F$. Let $P = MU$ be a proper parabolic subgroup of $G$ and let $\tau$ be an irreducible automorphic subrepresentation of $M(\BA)$. Let $\mathcal{A}^{\tau}_P(G)$ be the space consisting functions $\phi \in \mathcal{A}_P(G)$ such that for all $k \in K$ the function $m \mapsto e^{-\langle \rho_P, H_P(m)\rangle} \phi(mk)$ on $M(F)\backslash M(\BA)$ belongs to the space of $\tau$. Denote by $I_P^G(\tau)$ the representation space $\CA^{\tau}_P(G)$ equipped with the right translation action of $G(\BA)$. For $\phi \in \CA_P^{\tau}(G)$ and $\lambda \in \mathfrak{a}_{P,\BC}^{\ast}$, set
	\begin{align*}
		\phi_{\lambda} (g)  =  e^{\langle \lambda, H_P(g)\rangle} \phi(g), \quad g \in G(\BA).
	\end{align*}
	Denote by $I_P^G(\tau,\lambda)$ the representation space $\{\phi_{\lambda}\ |\ \phi \in \CA_P^{\tau}(G) \}$ equipped with  the right translation of  $G(\BA)$. For $P \subset Q$, the Eisenstein series $E^Q(\phi,\lambda)$ is defined by
	\begin{align}\label{formual::defn--Eisenstein}
		E^Q(g,\phi, \lambda) = \sum_{\gamma \in P(F) \backslash Q(F)} \phi_{\lambda} (\gamma g),\quad g \in G(\BA).
	\end{align}
	The series converges absolutely and uniformly in $g$ and $\lambda$ for $\Re\lambda$ sufficiently regular in the postive Weyl chamber of $\mathfrak{a}_P^{\ast}$. We still use $E^Q(\phi,\lambda)$ to denote its meromorphic continuation. When $Q = G$, we will write instead $E(\phi, \lambda)$ for simplicity. Note that $E^Q(\phi,\lambda) \in \CA_Q(G)$ and for all $k \in K$
	\begin{align*}
		E^Q(lk, \phi,\lambda) = E^L(l,R_k\phi,\lambda).
	\end{align*}
	Here $R_k\phi$ is the right translation of $\phi$ by $k$. We have $R_k \phi \in \CA_{L\cap P}^{\tau \otimes \chi}(L)$ with $\chi(m) = e^{\langle \rho_Q, H_0(m)\rangle}$.
	
	Denote by $\Omega = \Omega^G$ the Weyl group associated to $(G,A_0)$. For any two standard parabolic subgroups $P = M \ltimes U$ and $Q = L \ltimes V$, let 
	\begin{align*}
		{}_Q\Omega_P  = \{ w \in \Omega \ |\ w\alpha >0 \text{ for all }\alpha \in \Delta_0^P \text{ and } w^{-1}\alpha > 0 \text{ for all } \alpha \in \Delta_0^Q\}.
	\end{align*}
	Then for any $w \in {}_Q\Omega_P$, the group $M \cap w^{-1}Lw$ is the Levi subgorup of a tandard parabolic subgroup $P_w$ of $P$; the group $L \cap wMw^{-1}$ is the Levi subgroup of a standard parabolic subgroup $Q_w$ of $Q$. Set
	\begin{align*}
		&\Omega(P;Q) = \{ w \in {}_Q\Omega_P \ |\ Q_w = Q\},\\
		&\Omega(P,Q) = \{ w \in {}_Q\Omega_P \ |\ wMw^{-1} = L \} = \Omega(P;Q) \cap \Omega(Q;P)^{-1}.
	\end{align*}
	Note that if $w \in {}_Q\Omega_P$, then $w \in \Omega(P_w,Q_w)$.
	
	For any $w \in \Omega(P,Q)$, the intertwining operator $M(w,\lambda): \CA_P \ra \CA_Q$ is defined by the formula
	\begin{align*}
		(M(w,\lambda)\phi)_{w\lambda}(g) = \int_{(V \cap wUw^{-1})(\BA) \backslash V(\BA) } \phi_{\lambda} (w^{-1} u g) du, \quad g \in G(\BA).
	\end{align*}
	The integral converges locally uniformly in $g$ and $\lambda$ provided that $\langle \Re \lambda , \alpha^{\vee} \rangle \gg 0$ for every root $\alpha \in \Phi_P$ such that $w\alpha <0$. We still use $M(w,\lambda)$ to denote its meromorphic continuation. When $\lambda = 0$, we will write instead $M(w)$ for simplicity. Note that the operator $M(w,\lambda)$ maps $I_P^G(\tau,\lambda)$ into $I_Q^G(w\tau,w\lambda)$.
	
	The constant terms of Eisenstein series are computed as follows. We refer the readers to \cite[Section 5.10]{Bernstein-Lapid-meromorphic-continuation-Eisenstein} for the proof.
	\begin{lem}\label{lem::general constant term}
		We have
		\begin{align}\label{formula::general constant term}
			E(\phi,\lambda)_Q  = \sum_{w\in  {}_Q\Omega_P} E^Q(M(w,\lambda)\phi_{P_w},w\lambda), 			
		\end{align}
		where the summation appeared in the Eisenstein series on the right hand side is taken ove $Q_w(F) \backslash Q(F)$. For the cuspidal components of these constant terms, we have
		\begin{align}\label{formula::cusp-component-const-term}
			E^{\textup{cusp}}(\phi,\lambda)_Q  =\sum_{w \in \Omega(P;Q)} [M(w,\lambda)\phi_{P_w}]^{\textup{cusp}}_{w\lambda}= \sum_{w \in \Omega(P;Q)} [M(w,\lambda)(\phi^{\textup{cusp}}_{P_w})]_{w\lambda}.
		\end{align}
	\end{lem}
	
	\subsection{Constant terms}\label{sec::constant term Eisenstein series}
	
	The classification of the discrete spectrum for $\GL_n$ was established by Moeglin and Waldspurger \cite{Moeglin-Waldspurger-residue-spectrum-GLn}. We recall it here. Given a decomposition $n = dr$, let $P = M \ltimes U$ the standard parabolic subgroup associated to the partition $(r,\cdots,r)$. Set
	\begin{align*}
		\Lambda_d = ((d-1)/2, (d-3)/2,\cdots, (1-d)/2) \in (\mathfrak{a}_M^G)^{\ast}.
	\end{align*}
	Here we identity $\mathfrak{a}_{M,\BC}^{\ast}$ with the set of $d$-tuples. Let $\sigma$ be an irreducible cuspidal automorphic representation of $\GL_r(\BA)$ and $\tau = \sigma \otimes \cdots \otimes \sigma$ ($d$ times) be a cuspidal automorphic representation of $M(\BA)$. The induced representation $I_P^G(\tau,\Lambda_d)$ has a unique irreducible quotient that is denoted by $\speh(\sigma,d)$. The representation $\speh(\sigma,d)$ occurs in the discrete spectrum of $\GL_n$ with multiplicity one and can be realized using multi-residues of Eisenstein series as follows. For $\phi \in \CA_P^{\tau}(G)$, the meromorphic function
	\begin{align*}
		\prod_{i=1}^{d-1} (\lambda_i - \lambda_{i+1} - 1) E(\phi,\lambda)
	\end{align*}
	is holomorphic at $\lambda = \Lambda_d$. Denote by $E_{-1}(\phi)$ its limit when $\lambda \ra \Lambda_d$. This is a square-integrable automorphic form. The intertwining map $\phi_{\Lambda_d} \mapsto E_{-1}(\phi)$ factors through the quotient $I_P^G(\tau,\Lambda_d) \ra \speh(\sigma,d)$ and gives a realization of $\speh(\sigma,d)$ into the space of square-integrable automorphic forms of $\GL_n$. By convention, we still denote by $\speh(\sigma,d)$ this automorphic realization. Moreover, as we vary $r$ and $\sigma$ we get the entire discrete spectrum of $L^2([\GL_n])$. 
	
	Note that the set ${}_P\Omega_P$ can be identified naturally with the group $\mathfrak{S}$ of permutations on $d$ elements. Let $Q$ be the parabolic subgroup associated to the partition $(m_1r,\cdots,m_kr)$ with $m_1 + \cdots + m_k = d$. For $\phi \in \CA^{\tau}_P(G)$, the singularities of $E^Q(\phi,\lambda)$ are the singularities of its cuspidal components. By \eqref{formula::cusp-component-const-term}, they are contained in the singularities of $M(w,\lambda)$ where $w \in \Omega(P,P) \cap \Omega^Q$, which are along the hyperplanes 
	\begin{align*}
		\lambda_i - \lambda_j - 1 = 0, \quad  i,j \text{ are in the same segment and }  i < j.
	\end{align*}
	Here, by saying $i,j$ are in the same segment, we mean that there exists $1\leqslant l \leqslant k$ such that $$i,j \in  \left[\sum_{i=1}^{l-1} m_i+1, \sum_{i=1}^{l}m_i \right].$$
	Set
	\begin{align*}
		\Lambda^Q &  =  ( \Lambda_{m_1},\cdots ,\Lambda_{m_k}) \in (\mathfrak{a}_P^Q)^{\ast} .
	\end{align*}
	Let $\mu \in \mathfrak{a}_Q^{\ast}$, we define
	\begin{align*}
		E_{-1}^Q(\phi,\mu) = \lim_{\lambda \ra \Lambda^Q} \prod_{i \in \mathcal{I}^Q}(\lambda_i - \lambda_{i+1} -1) E^Q(\phi, \lambda + \mu),
	\end{align*}
	where $\mathcal{I}^Q$ stands for the set $\{ 1,2,\cdots,d-1\} \backslash \{ m_1,m_1+m_2,\cdots,\sum_1^{k-1} m_l\}$. It defines a surjective intertwining operator 
	\begin{align*}
		E^Q_{-1}(\mu) \colon I_P^G(\tau,\Lambda^Q + \mu) \ra I_Q^G(\speh(\sigma,m_1)\otimes \cdots \otimes \speh(\sigma,m_k),\mu).
	\end{align*}
	
	We also define the multi-residue of the intertwining operator $M(w,\lambda)$ as follows. Set
	\begin{align*}
		M_{-1}(w) = \lim_{\lambda \ra \Lambda_d} 	M (w,\lambda) \prod_{1\leqslant i < d,\  w(i) > w(i+1)} (\lambda_i - \lambda_{i+1} -1) .
	\end{align*}
	It is an intertwining operator from $I_P^G(\tau,\Lambda_d)$ to $I_P^G(\tau,w\Lambda_d)$.
	
	Define by $w_Q$ the longest element in $\Omega(Q;P)^{-1} = \{ w \in {}_Q\Omega_P \ |\ P_w = P\}$, which can be identified with the permutation
	\begin{align*}
			\begin{pmatrix}
					1 &\cdots & m_k   \cdots & d-m_1+1 &\cdots & d  \\
					d-m_k+1 & \cdots & d   \cdots & 1 & \cdots & m_1
				\end{pmatrix}.
		\end{align*}

	\begin{lem}\label{lem::constant term--residue}
		Let notations be as above. Let $\phi \in \CA_P^{\tau}(G)$ such that $E_{-1}(\phi) \in \speh(\sigma,d)$. Then $E_{-1}(\phi)_Q =0$ unless $P \subset Q$. When $P \subset Q$, we have
		\begin{align*}
			E_{-1}(\phi)_Q  = E^Q_{-1}(M_{-1}(w_Q)\phi, -(\Lambda_d)_Q)
		\end{align*}
	\end{lem}
	\begin{proof}
		The special case where $k=2$ have been proved in \cite[Lemma 2.4]{Offen-Sayag-Global-Klyachko} and \cite[Lemma 6.1]{Yamana-Residue-Periods-JFA}. By \eqref{formula::general constant term}, the constant term of the Eisenstein series $E(\phi,\lambda)$ along $Q$ is given by
		\begin{align*}
			E(\phi,\lambda)_Q = \sum_{w \in \Omega(Q;P)^{-1}} E^Q (M(w,\lambda)\phi, w\lambda).
		\end{align*}
		Clearly, if $P$ is not contained in $Q$,  we have $E_{-1}(\phi)_Q = 0$. For $P \subset Q$, we need to show that after applying the multi-residue operator, only the term associated to $w_Q$ is nonzero. Assume $Q$ is of type $(m_1r,\cdots,m_kr)$. For each $w \in \Omega(Q;P)^{-1}$, set
		\begin{align*}
			\mathcal{I}_w^{(1)} &= \{ i \in [1,d-1] \ |\ wi > w(i+1) \}  \\
			\mathcal{I}_w^{(2)} &= \{ i\in [1,d-1] \ |\  wi < w(i+1)\text{, $wi$ and $w(i+1)$ are in the same segment}  \}.
		\end{align*}
		Note that, by definition,
		\begin{align*}
			M_{-1}(w)\phi = \lim_{\lambda \ra \Lambda_d} \prod_{i \in \mathcal{I}_w^{(1)}} (\lambda_i - \lambda_{i+1} - 1) M(w,\lambda) \phi
		\end{align*}  
		Also note that, 
		\begin{align*}
			\prod_{i \in \mathcal{I}_w^{(2)}}(\lambda_i - \lambda_{i+1} - 1)E^Q(M_{-1}(w)\phi,w\lambda)
		\end{align*}
		is holomorphic at $\lambda = \Lambda_d$. We show that if $\mathcal{I}_w^{(1)} \cup \mathcal{I}_w^{(2)} = [1,d-1]$, then $w = w_Q$. As $w \in {}_Q\Omega_P$ , we have $w^{-1}(1) < w^{-1}(2) < \cdots < w^{-1}(m_1)$. Note that $w^{-1}(m_1)$ lies in neither $\mathcal{I}_w^{(1)}$ nor $\mathcal{I}_w^{(2)}$, so $w^{-1}(m_1) = d$. The above inequality also implies that all $w^{-1}(j)$, $1\leqslant j \leqslant m_1 -1$, lie in $\mathcal{I}_w^{(2)}$. This forces that $w^{-1}(1) = d - m_1 + 1$ and $w^{-1}(j+1) = w^{-1}(j) + 1$ for $ 1 \leqslant j \leqslant m_1 -1$. By repeating the above argument, we finally find that $w = w_Q$. Hence
		\begin{align*}
			E_{-1}(\phi)_Q = \lim_{\lambda \ra \Lambda_d} \prod_{i \in \mathcal{I}^Q}(\lambda_i - \lambda_{i+1} - 1)E^Q(M_{-1}(w_Q)\phi,w_Q\Lambda_d).
		\end{align*}
		So we get the equality in the lemma after noting that
		\begin{align*}
			w_Q\Lambda_d  - \Lambda^Q = -(\Lambda_d)_Q.
		\end{align*}
	\end{proof}
	
	We end this section by computing the constant terms of the mirabolic Eisenstein series $\tilde{E}(\Phi,s)$, $\Phi \in \mathcal{S}(\BA^n)$. Let $P=MU$ be a standard parabolic subgroup of $G_n$. By Lemma \ref{lem::general constant term}, we can write
	\begin{align}\label{formual::constant term--mirabolic Eisenstein}
		\tilde{E}_P(\Phi,s) = \sum_{w\in {}_P\Omega_{\mathcal{P}_n}} \tilde{E}^P(M(w)F(\Phi;s)).
	\end{align}
	We identify the Weyl group $\Omega$ of $G_n$ with the permutation group $\mathfrak{S}_n$ and denote by $w_i$ the cyclic permutation $(i,i+1,\cdots,n)$. Assume $P$ is of  type $(n_1,\cdots,n_t)$. Then ${}_PW_{\mathcal{P}_n}$ can be identified with the set $\{ w_i \ |\ i \in \mathcal{I}_P \}$ with $\mathcal{I}_P = \{ n_1 + \cdots + n_j \ |\ 1\leqslant j \leqslant t\}$, or written differently, with the set $\{w_{(i)} \ |\ i =1,\cdots,t \}$, where $w_{(i)}$ is the cyclic permutation
	\begin{align*}
		(\sum_1^i n_j, \sum_1^i n_j + 1 ,\cdots,n).
	\end{align*}
	
	\begin{lem}\label{lem::const-term-mirabolic}
		Let $m = \textup{diag}(m_1,\cdots,m_t) \in M$. Then
		\begin{align*}
			\tilde{E}^P(mg,M(w_{(i)})F(\Phi;s))  =   e^{\langle \rho_0 + w_{(i)}(2s \rho_{\mathcal{P}_n} - \rho_0),H_P(m) \rangle} | g|^s \tilde{E}(m_i,\Phi';s') 
		\end{align*}
		for $g \in G_n(\BA)$. Here $s' = (ns - \sum_{j=i+1}^t n_j)/n_i$ and $\Phi' \in \CS(\BA^{n_i})$ is given by 
		\begin{align*}
			\Phi'(z) = \int_{\BA^{n_{i+1}+ \cdots + n_t}} \Phi((\underbrace{0,\cdots,0}_{n_1 + \cdots + n_{i-1}},z,x)g) dx.
		\end{align*}
	\end{lem}
	\begin{proof}
		It is known that
		\begin{align*}
			\tilde{E}^P(e^Xg,M(w_{(i)})F(\Phi;s)) = e^{\langle \rho_0 + w_{(i)}(2s \rho_{\mathcal{P}_n} - \rho_0), X \rangle} \tilde{E}^P(g,M(w_{(i)})F(\Phi;s))
		\end{align*}
		for $X \in \mathfrak{a}_P$ and $g \in \GL_n(\BA)$ \cite[Lemma 6.3]{Yamana-Residue-Periods-JFA}. Hence we may assume that $H_P(m) = 0$. Note that $P_{w_{(i)}}$ is the standard parabolic subgroup associated to the composition $(n_1,\cdots,n_{i-1},n_i - 1, 1,\cdots,n_t)$. Thus $\tilde{E}^P(mg,M(w_{(i)})F(\Phi;s))$ is the sum over $\gamma \in P_{(n_i - 1,1)}(F) \backslash \GL_{n_i}(F)$ of 
		\begin{align*}
			&\int_{(U_{P_{w_{(i)}}} \cap w_{(i)} U_{(\mathcal{P}_n)_{w_{(i)}}} w_{(i)}^{-1} )(\BA) \backslash U_{P_{w_{(i)}}}(\BA)} F(w_{(i)}^{-1} u \gamma mg,\Phi;s) du \\
			= & | g|^s \int_{\BA^{n_{i+1}}} \cdots \int_{\BA^{n_t}} \int_{\BA^{\times}} \Phi(a (\underbrace{0,\cdots,1}_{n_1 + \cdots + n_i},\, x_{i+1},\cdots,x_t) \gamma mg)|a|^{ns} d^{\times} a dx_{i+1}\cdots dx_t  \\
			= & |  g|^s \int_{\BA^{n_{i+1}+ \cdots + n_t}}\int_{\BA^{\times}}\Phi( (\underbrace{0,\cdots,0}_{n_1+\cdots +n_{i-1}}, ae_{n_i}\gamma m_i,x)g) |a|^{ns-\sum_{j=i+1}^t n_j}  d^{\times} a dx \\
			= & |  g|^s F(\gamma m_i,\Phi';s').
		\end{align*}
	\end{proof}
	

	\section{Global linear periods: The square-integrable case}\label{sec::linear period-square-integral case}
	
	In this section, let $G = G_{2n}$ and $G' = G'_{2n}$. We retain the notations of previous sections. As the quotient $G'(F)\backslash G'(\BA)^{1,G}$ has infinite volume, the period integral \eqref{formula::defn-period-integral} may diverge for square-integrable automorphic forms. This is obvious for constant functions. We continue to use $\period(\varphi,s)$ for the regularized linear period of an automorphic form $\varphi$, in the sense of Zydor, whenever it is defined. In this case we also call $\varphi$ is $\mu_s$ regular.
	
	We explicate now the set $\CF^{G}(P_0')$ of relative parabolic subgroups. Using dynamic theory of parabolic subgroups \cite[Section 2.4]{Zydor-periods-AENS}, it is not hard to check that $\CF^G(P_0')$ is in bijection with the set of datum
	\begin{align*}
		(n_1,\cdots,n_k;n_1',\cdots,n_k'),\quad 1 \leqslant k \leqslant 2n,
	\end{align*}
	satisfying
	\begin{align*}
		n_i,\ n_i' \geqslant 0,\ n_i + n_i' > 0 \text{ for all $1 \leqslant i \leqslant k$}
	\end{align*}
	and 
	\begin{align*}
		n = \sum_{j=1}^{k} n_j = \sum_{j=1}^{k} n_j'.
	\end{align*}
	If $P \in \CF^G(P_0')$ correspondes to $(n_1,\cdots,n_k;n_1',\cdots,n_k')$, then $P$ is conjugate to the standard parabolic subgroup $P_{(n_1+n_1',\cdots,n_k+n_k')}$. We choose and will fix $w_P \in \Omega^G$ such that $w_PPw_P^{-1} = P_{(n_1+n_1',\cdots,n_k + n_k')}$.
	
	Suppose that $2n = dr$ and $\sigma$ is an irreducible cuspidal automorphic representation of $G_r(\BA)$ with trivial central character.
	 
	\begin{lem}\label{lem::phi*E--regular}
		Let $\varphi \in \speh(\sigma,d)$. The pair $(\varphi,E(\Phi,s_1,s))$ belongs to $\CA(G \times G')^{\ast\ast}$ for $(s_1,s) \in \BC^2$ outside a finite union of hyperplanes. 
	\end{lem}
	\begin{proof}
		We need to show that, for $P \subset R \subset Q$ in $\CF^G(P_0')$ with $\dim \mathfrak{a}_R^Q = 1$,
		\begin{align}\label{formula::condition--xx}
			\langle \lambda  + \lambda' + \rho_{P'} + \underline{\rho}_P, \mathfrak{a}_R^Q \rangle \ne 0
		\end{align}
		for all $\lambda \in \mathcal{E}_P(\varphi)$ and $\lambda' \in \mathcal{E}_{P'}(E(\Phi;s_1,s))$.
		
		Suppose that 
		\begin{align*}
			P' = w_{\star}  \begin{pmatrix}
				P_1 &  \\
				    & P_2
			\end{pmatrix}  w_{\star}^{-1},
		\end{align*} 
		where $P_1$ and $P_2$ are two standard parabolic subgroups of $G_n$. By \eqref{formula::connection-Eisenstein-Mirabolic}, we have
		\begin{align*}
			E_{P'}(\iota(g_1,g_2),\Phi;s_1,s) = | g_1|^{s_1 + s - 1/2} | g_2|^{-s_1 - s + 1/2} \tilde{E}_{P_2}(g_2,\Phi,2s_1).
		\end{align*}
		By Lemma \ref{lem::const-term-mirabolic}, $\lambda' + \rho_{P'}$, $\lambda' \in \mathcal{E}_{P'}(E(\Phi;s_1,s))$, are of the form
		\begin{align}\label{formula::exponent-Eisenstein}
			w_{\star}((s_1 + s -\frac{1}{2})(\underbrace{1,\cdots,1}_n,\underbrace{-1,\cdots,-1}_n)+(\underbrace{0,\cdots,0}_n,\,\underbrace{\rho_{0,n}+w_{i,n}(4s_1\rho_{\mathcal{P}_n}-\rho_{0,n})}_n)),
		\end{align}
		with $i \in \mathcal{I}_{P_2}$ and $w_{i,n} \in \Omega^{G_n}$ the cyclic permutation introduced in the last section. Note that 
		\begin{align*}
			4w_{i,n}\rho_{\mathcal{P}_n} = (\underbrace{2,\cdots,2,2-2n}_i,2,\cdots,2)
		\end{align*}
		and 
		\begin{align*}
			\rho_{0,n} - w_{i,n}\rho_{0,n} = (\underbrace{0,\cdots,0,n-i}_i,-1,\cdots,-1).
		\end{align*}
		
		Suppose that $R$ is conjugate to the standard parabolic subgroup $P_{(n_1,n_2,\cdots,n_k)}$ and $Q$ is conjugate to $P_{(n_1,\cdots,n_j + n_{j+1},\cdots,n_k)}$ for some $1 \leqslant j \leqslant k-1$. Then $\mathfrak{a}_R^Q$ is generated by 
		\begin{align}\label{formula::defn-t_P^Q}
			t_R^Q := w_R^{-1} (\underbrace{0,\cdots,0}_{n_1 + n_2 +\cdots + n_{j-1}},\underbrace{n_{j+1},\cdots,n_{j+1}}_{n_j},\underbrace{-n_j,\cdots,-n_j}_{n_{j+1}},0,\cdots,0).
		\end{align}
		Suppose that $R$ correspondes to $(n_1',\cdots,n_k';n_1'',\cdots,n_k'')$, then
		\begin{align*}
			\langle w_{\star}(\underbrace{1,\cdots,1}_n,\underbrace{-1,\cdots,-1}_n), t_R^Q \rangle = 2(n_j'n_{j+1}'' - n_{j+1}'n_j'').
		\end{align*}
		If $n_j'n_{j+1}'' - n_{j+1}'n_j'' \ne 0$, then the condition \eqref{formula::condition--xx} are satisfied for $(s_1,s)$ outside a finite union of hyperplanes. Now we assume that $n_j'n_{j+1}'' - n_{j+1}'n_j'' = 0$. Note that $n_{j+1}'$ and $n_{j+1}''$ are not all zero, so there exist a positive rational number $c$ such that $n_j' = c n_{j+1}'$ and $n_j'' = c n_{j+1}''$. Note that 
		\begin{align*}
			(\underbrace{1,\cdots,1}_n,\underbrace{-1,\cdots,-1}_n) + (\underbrace{0,\cdots,0}_n,\,\underbrace{ 4 w_{i,n} \rho_{\mathcal{P}_n}}_n )
		\end{align*}
		equals to
		\begin{align*}
			(\underbrace{1,\cdots,1}_n,\underbrace{1,\cdots,1,1-2n}_{i},1,\cdots,1),
		\end{align*}
		for $i \in \mathcal{I}_{P_2}$. If 
		\begin{align*}
			\langle w_{\star}((\underbrace{1,\cdots,1}_n,\underbrace{1,\cdots,1,1-2n}_{i},1,\cdots,1)), t_R^Q \rangle \ne 0,
		\end{align*}
		then the condition \eqref{formula::condition--xx} are satisfied for $(s_1,s)$ outside a finite union of hyperplanes. If, otherwise, 
		\begin{align*}
			\langle w_{\star}((\underbrace{1,\cdots,1}_n,\underbrace{1,\cdots,1,1-2n}_{i},1,\cdots,1)), t_R^Q \rangle = 0,
		\end{align*}
		then
		\begin{align*}
			\langle w_{\star}((\underbrace{0,\cdots,0}_n,\,\underbrace{\rho_{0,n} - w_{i,n}\rho_{0,n}}_n)), t_R^Q \rangle = 0.
		\end{align*}
		Since $\varphi$ has a single cuspidal exponent $-\Lambda_d$ relative to standard parabolic subgroups. It remains to show that
		\begin{align}\label{formula::twist-regular}
			\langle - w_P^{-1} \Lambda_d + \underline{\rho}_P, t_R^Q \rangle \ne 0.
		\end{align}
		Note that the projection of $\underline{\rho}_P$ to $\mathfrak{a}_R^Q$ is $\underline{\rho}_R$. We have
		\begin{align*}
			\langle \underline{\rho}_R, t_R^Q \rangle = (-c^2 - c) n_{i+1} (n_{i+1}' - n_{i+1}'')^2/2 \leqslant 0.
		\end{align*}
		It is well known that $ \langle - w_P^{-1}\Lambda_d, t_R^Q \rangle = \langle - \Lambda_d, w_R t_R^Q \rangle < 0$. So we are done.
	\end{proof}

	\begin{rem}
		If we let $s = 0$, it is possible that $\varphi \otimes E(\Phi;s_1,0)$ does not belong to $\CA(G \times G')^{\ast\ast}$ for any $s_1 \in \BC$. For example, let $G = G_{18}$, $r = 6$, $P$ correspondes to $(5,3,1;1,3,5)$ and $Q$ correspondes to $(8,1;4,5)$. The condition \eqref{formula::condition--xx} is not satisfied. Explicitly, the term 
		\begin{align*}
			\int_{Q'(F) \backslash G'(\BA)^{1,G}}^* \varphi_Q (g) E_{Q'}(g;\Phi;s_1,0) \hat{\tau}_P (H_{0'}(g)^{G}_{P} - T^G_P) dg.
		\end{align*}
		is not defined.	This is one of the reasons why we have to introduce the parameter $s$.
	\end{rem}
	
	\begin{lem}\label{lem::regularizable--analytic family}
		Fix $\Phi \in \CS(\BA^n)$. Let $\varphi(\lambda)$ be an analytic family of automorphic forms (not necessarily square-integrable). Let $\mathcal{O}$ denote the set of all triplets $(\lambda,s_1,s)$ such that
		\begin{align*}
			(\varphi(\lambda),E(\Phi;s_1,s)) \in \CA(G  \times G' )^{\ast}.
		\end{align*}
		Then $\mathcal{O}$ is nonempty and open and $(\lambda,s_1,s) \mapsto \period^{G'}(\varphi(\lambda)\otimes E(\Phi;s_1,s))$ is an analytic function on $\mathcal{O}$.
	\end{lem}
	\begin{proof}
		For a fixed $\lambda$, we have to check the condition \eqref{formula::condition--xx} with $Q = G$. Note that
		\begin{align*}
			\langle w_{\star}((\underbrace{1,\cdots,1}_n,\underbrace{1,\cdots,1,1-2n}_{i},1,\cdots,1)), t_R^G \rangle \ne 0,
		\end{align*}
		for all $R \in \CF^G(P'_0)$ with $ \dim \mathfrak{a}_R^G = 1$. Hence $(\varphi(\lambda),E(\Phi;s_1,s))$ belongs to $\CA(G \times G')^{\ast}$ except for a finite number of $s_1$. By \cite[Theorem 3.9]{Zydor-periods-AENS}, the integral $\period^{G',T}(\varphi(\lambda) \otimes E(\Phi;s_1,s))$ is uniformly convergent for $\lambda$, $s_1$ and $s$ in compact subset, which completes the proof.
	\end{proof}
	
	Actually we have the following vanishing result.
	\begin{lem}\label{lem::vanishing--phi*E}
		If $d \geqslant 2$, then $\period^{G'}(\varphi \otimes E(\Phi;s_1,s))$ is identically zero for all $\varphi \in \text{Sp}(\sigma,d)$ and $\Phi \in \CS(\BA^n)$.
	\end{lem}
	\begin{proof}
		We still write $R$ for the right translation of $G'(\BA)$ on $\mathcal{S}(\BA^n)$ through the second factor. Note that $\Phi \mapsto f(\Phi;s_1,s)$ is an $G'(\BA)$-intertwining map. By Lemma \ref{lem::regularizable--analytic family} and part (ii) of Proposition \ref{prop::period-well defined}, the composition of this map with the regularized period defines an element in the space $\Hom_{G'(\BA)} (\speh(\sigma,d) \otimes R, \nu^a \boxtimes \nu^b)$ for some $a,b$ which are affine functions in $s_1$ and $s$. Locally such invariant functional does not exist by Theorem \ref{thm::generic uniqueness} for $s_1,s$ in a general position. Thus the meromorphic function $\period^{G'}(\varphi \otimes E(\Phi;s_1,s))$ vanishes identically.
	\end{proof}

	We define functions on $G'(\BA)$ by
	\begin{align*}
		W^{\psi}_r (g,\varphi)  &= \int_{N_r(F) \backslash N_r(\BA)} \varphi_{P_{(2n-r,r)}} \left( \begin{pmatrix}
			\mathbf{1}_{2n-r}  &  \\  & u 
		\end{pmatrix}  g \right) \overline{\psi(u)} du, \\
		\mathrm{\Theta}^{\psi}(g,\varphi;s) &= \int_{G'_{2n-r}(F) \backslash G'_{2n-r}(\BA)^{1,G_{2n-r}}}^{\ast} W_r^{\psi} \left( \begin{pmatrix}
			m &  \\  & \mathbf{1}_r
		\end{pmatrix} g,\varphi\right) \mu_s(m) dm
	\end{align*} 
	
	For a semi-standard parabolic subgroup $P$, we denote by $v_P$ the volume of the parallelogram
	\begin{align}\label{formula::defn-v_P-volume}
		\left\{   \sum_{\alpha \in \Delta_P} a_{\alpha} \alpha \ |\ 0 \leqslant a_{\alpha} \leqslant 1  \right\}.
	\end{align}
	in $\mathfrak{a}_P^G$.
	
	\begin{prop}\label{prop::period formula}
		Let $\varphi \in \speh(\sigma,d)$, $d \geqslant 2$. Then for all but finitely many $s \in \BC$, $\varphi$ is $\mu_s$-regular. $\period(\varphi,s)$ can be analytically extended to all of $\BC$. Moreover, for any $\Phi \in \CS(\BA^n)$, we have
		\begin{align*}
			\hat{\Phi}(0)\period(s,\varphi) = r v_P \int_{(P^{(r)}_{2n})'(\BA) \backslash G_{2n}'(\BA)}  \mathrm{\Theta}^{\psi}(g,\varphi ; s-1/2) \Phi(e_n g_2) \mu_s(g) |\det g_2|dg.
		\end{align*}
		Here $P$ is the standard parabolic subgroup of type $(2n-r,r)$, and we write $g = w_{\star} \textup{diag}(g_1,g_2) w_{\star}^{-1}$.
	\end{prop}
	\begin{proof}
		The first assertion can be proved by \eqref{formula::twist-regular}. For $T$ sufficiently positive, by Lemma \ref{lem::pole--mirabolic Eisenstein}, 
		\begin{align}\label{formula::twisted truncated period as residue}
			\lim_{s_1 \ra 1/2} (s_1 - \frac{1}{2}) \int_{[G']^{1,G}} \Lambda^T \varphi (g) E(g,\Phi;s_1,s) dg =v_F \frac{\hat{\Phi}(0)}{2n}\int_{[G']^{1,G}} \Lambda^T\varphi(g)\mu_{s}(g) dg .
		\end{align}
		By Lemma \ref{lem::phi*E--regular} and Proposition \ref{prop::TruncatedPeriod-Main}, it is permissible to write 
		\begin{align*}
			\int_{[G']^{1,G}} \Lambda^T \varphi (g) E(g,\Phi;s_1,s) dg
		\end{align*}
		as the sum over $P \in \CF^G(P_0')$ of 
		\begin{align*}
			J_P =  \varepsilon_P^G \int_{P'(F) \backslash G'(\BA)^{1,G}}^{\ast} \varphi_P (g) E_{P'}(g,\Phi;s_1,s) \hat{\tau}_P(H_{0'}(g)_P^G - T_P^G) dg.
		\end{align*}
		By Lemma \ref{lem::vanishing--phi*E}, $J_G = 0$. Assume that $P$ correspondes to $(n_1',\cdots,n_k';n_1'',\cdots,n_k'')$ and that \begin{align*}
			P' = w_{\star}\begin{pmatrix}
				P_1 &  \\
				& P_2
			\end{pmatrix}w_{\star}^{-1}
		\end{align*}
		with $P_1,\ P_2$ two standard parabolic subgroups of $G_n$.  We also assume that $n'_i + n_i'' = m_ir$ for $i = 1,\cdots,k$, otherwise $\varphi_P = 0$ by Lemma \ref{lem::constant term--residue}. By the formula \eqref{formula::connection-Eisenstein-Mirabolic} and \eqref{formual::constant term--mirabolic Eisenstein}, 
		\begin{align*}
			E_{P'}((g_1;g_2),\Phi;s_1,s) &= |g_1|^{s_1 + s- 1/2} | g_2 |^{-s_1 -s + 1/2} \tilde{E}_{P_2}(g_2,\Phi;2s_1)  \\
			&= \sum_{i \in \mathcal{I}_{P_2}} |g_1|^{s_1 + s- 1/2} |g_2 |^{-s_1 -s + 1/2} \tilde{E}^{P_2} (g_2,M(w_{i,n})F(\Phi;2s_1))  \\
			&:= \sum_{i \in \mathcal{I}_{P_2}} E_{P',i}((g_1;g_2),\Phi;s_1,s).
		\end{align*}
		Then 
		\begin{align}\label{formula::J_P}
			\begin{aligned}
				J_P = \varepsilon_P^G \sum_{i \in \mathcal{I}_{P_2}} & \int_{K'} \int_{[M']^{1,M}}^* \varphi_P(mk) E_{P',i}(mk,\Phi;s_1,s) e^{- \langle 2\rho_{P'}, H_{P'}(m) \rangle} dm  dk \\
				& \cdot \int_{\mathfrak{a}_P^G}^{\sharp}  e^{\langle -w_P^{-1}\Lambda_d + \lambda'_i(s_1,s) + \underline{\rho}_P, X \rangle} \hat{\tau}_P(X - T_P^G) dX,
			\end{aligned}
		\end{align}
		where $\lambda'_i(s_1,s)$ stands for the exponent given in \eqref{formula::exponent-Eisenstein}, that is,
		\begin{align*}
			w_{\star}((s_1 + s -\frac{1}{2})(\underbrace{1,\cdots,1}_n,\underbrace{-1,\cdots,-1}_n)+(\underbrace{0,\cdots,0}_n,\,\underbrace{\rho_{0,n}+w_{i,n}(4s_1\rho_{\mathcal{P}_n}-\rho_{0,n})}_n)).
		\end{align*}
		The $\sharp$-integral in $J_P$ is equal to
		\begin{align}\label{formula::J_P-sharp}
			\frac{\varepsilon_P^G \, v_P\,e^{\langle -w_P^{-1}\Lambda_d + \lambda'_i(s_1,s) + \underline{\rho}_P, T_P^G \rangle} }{\prod_{\alpha \in \Delta_P} \langle -w_P^{-1}\Lambda_d + \lambda'_i(s_1,s) + \underline{\rho}_P, \alpha \rangle},
		\end{align}	
		by the explicit formula in \cite[(15)]{J-L-R-periods-JAMS}. Note that each term in $J_P$ is meromorphic in $s_1$. After taking the residue at $s_1 = 1/2$, the left hand side of \eqref{formula::twisted truncated period as residue} is a polynomial exponential function in $T$ with exponents contained in the set $\{ (-w_P^{-1}\Lambda_d + \lambda'_i(1/2,s) + \underline{\rho}_P)_P^G\}$. If 
		  \begin{align*}
		  	w_{\star}((\underbrace{1,\cdots,1}_n,\underbrace{-1,\cdots,-1}_n))
		  \end{align*}
		 is not zero on $\mathfrak{a}_P^G$, then for all but finitely many $s \in \BC$, the exponents are nonzero. Hence the $J_P$ term do not contribute to the purely polynomial part of the left hand side of \eqref{formula::twisted truncated period as residue}. Assume that 
		\begin{align*}
			\langle w_{\star}((\underbrace{1,\cdots,1}_n,\underbrace{-1,\cdots,-1}_n)),\mathfrak{a}_P^G \rangle = 0.
		\end{align*}
		This is equivalent to the condition $n_i' = n_i''$ for $i = 1, \cdots,k$, which implies that $P$ is a standard parabolic subgroup. Note that
		\begin{align*}
			\langle -\Lambda_d + \underline{\rho}_P,\, (\underbrace{1,\cdots,1}_{m_1r},0,\cdots,0) \rangle = - m_1 (2n-m_1 r)/2 < 0
		\end{align*}
		and 
		\begin{align*}
			\langle \lambda_i'(1/2,s),\, (\underbrace{1,\cdots,1}_{m_1r},0,\cdots,0) \rangle = 0 \text{ or }[\frac{m_1r}{2}],
		\end{align*}
		depending on $m_1r \geqslant 2i$ or $m_1r < 2i$. Thus, for the purely polynomial part of the left hand side of \eqref{formula::twisted truncated period as residue} is nonzero, it is necessary that $P$ is the standard parabolic subgroup $P_{(2n-r,r)}$, $r$ is an even positive integer and $i = n \in \mathcal{I}_{P_2}$.
		
		We have shown that the purely polynomial part of the right hand side of \eqref{formula::twisted truncated period as residue} occur only when $P = P_{(2n-r,r)}$ and $i = n$. Note that in this case the regularized integral over $[M']^{1,M}$ in \eqref{formula::J_P} is holomorphic at $s_1 = 1/2$. This fact can be verified by arguments similar to those in Lemma \ref{lem::regularizable--analytic family}. Note that $\Delta_P = \{ \alpha \}$ with 
		\begin{align*}
			\alpha = \Big( \underbrace{\frac{1}{2n-r},\cdots,\frac{1}{2n-r}}_{2n-r},\underbrace{-\frac{1}{r},\cdots,-\frac{1}{r}}_r \Big).
		\end{align*}
		The residue of \eqref{formula::J_P-sharp} with $i = n$ at $s_1 = 1/2$ is $1/d$. Note also that, for $m \in [M']^{1,M}$,
		\begin{align*}
			\langle 2\rho_{P'},H_{P'}(m) \rangle = \langle \rho_P,H_{P'}(m) \rangle = 0.
		\end{align*}
		Thus, we have
		\begin{align}\label{formula::period = Levi P part}
			\hat{\Phi}(0) \period(\varphi,s) = r \frac{v_P}{v_F}\int_{K'} \int_{[M']^{1,M}}^* \varphi_P(mk) E_{P',n}(mk,\Phi;1/2,s)  dm  dk.
		\end{align}
		For $k \in K'$, write $k = w_{\star}^{-1} \textup{diag}(k_1,k_2) w_{\star}^{-1}$. For $m \in [M']^{1,M}$, write 
		\begin{align*}
			 m = w_{\star}  \left( \begin{matrix}
			 	m_1 & & & \\
			 	 & h_1 & & \\
			 	  & & m_2 & \\
			 	  & & & h_2
			 \end{matrix}\right) w_{\star}^{-1},
		\end{align*}
	    where $m_1,m_2 \in G_{n-r/2}(\BA)$ and $h_1,h_2 \in G_{r/2}(\BA)$. Note that we have $|m_1m_2| = |h_1 h_2| = 1$. Then, by \eqref{formula::connection-Eisenstein-Mirabolic} and Lemma \ref{lem::const-term-mirabolic},
		\begin{align*}
			E_{P',n}(mk,\Phi;1/2,s) & =  | m_1 h_1 |^s |  m_2|^{-s+1} |  h_2 |^{-s + 1 -d} \tilde{E}(h_2,\Phi_1;d) \\
			& = |m_1|^{s-1/2} |m_2|^{-s+1/2} E ((h_1;h_2),\Phi_1;d/2,s)
		\end{align*}
		where $\Phi_1\in S(\BA^{r/2})$ is given by 
		\begin{align*}
			\Phi_1(z) = \Phi((0,z)k_2).
		\end{align*}
		Define $\phi \in \mathcal{A}_P(G)$ by
		\begin{align}\label{formula::phi--varphi_P}
			\phi (g) = e^{\langle -\rho_P + (\Lambda_d)_P,\, H_P(g) \rangle} \varphi_P(g).
		\end{align}
		Hence, the inner integral over $[M']^{1,M}$ in the right hand side of \eqref{formula::period = Levi P part} is equal to
		\begin{align}\label{formula::period-1st}
			&\int^{\ast}_{[\GL_{2n-r}']^{1,\GL_{2n-r}}} \left( \int_{[\GL_r']^{1,\GL_r}} \phi \left( \begin{pmatrix}
				m & \\
				  & h 
			\end{pmatrix} k \right) E(h,\Phi_1;d/2,s) dh  \right)  \mu_{s-1/2}(m) dm.
		\end{align} 
		By Lemma \ref{lem::constant term--residue}, for any $k \in K'$, the function $\textup{diag}(m,h) \mapsto \phi(\textup{diag}(m,h)k) $ belongs to $\speh(\sigma,d-1) \otimes \sigma$. We should be careful as both sides of this equality include regularized integrals. We can verify this ``Fubini-like'' identity as follows. Replace the function $\phi(\cdot k)$ on $M$ by its trucation, use Fubini's theorem and then take the purely polynomial part. The inner integral of \eqref{formula::period-1st} converges absolutely and gives an analytic family of functions in $\speh(\sigma,d-1)$. By induction arguments, we could conclude that $\period(\varphi,s)$ can be analytically extended to $\BC$. By Proposition \ref{prop::Bump-F}, the inner integral of   \eqref{formula::period-1st} is equal to
		\begin{align*}
			v_F \int_{N_{r/2}(\BA) \backslash G_{r/2}(\BA)}  \int_{N_{r/2}(\BA) \backslash G_{r/2}(\BA)} W_r^{\psi} &\left( \begin{pmatrix}
				m &  \\
				  & h 
			\end{pmatrix} k,\phi\right)  \Phi((0,e_{r/2}h_2)k_2) \\ 
		 &\quad \cdot | h_1 |^{d/2+s -1/2} | h_2 |^{d/2-s+1/2} d h_1 d h_2 .
		\end{align*}
		Hence, $ \hat{\Phi}(0) \period(\varphi,s)$ is equal to $rv_P$ times 
		\begin{align}\label{formula::period-2nd}
			\int_{K'} \int_{(N_{r/2}(\BA) \backslash G_{r/2}(\BA))^2} &\left( \int^{\ast}_{[G_{2n-r}']^{1,G_{2n-r}}}W_r^{\psi} \left( \begin{pmatrix}
				m &  \\
				& h 
			\end{pmatrix} k,\phi\right) \mu_{s-1/2}(m) d m \right)  \\
		&	\cdot  \Phi((0,e_{r/2}h_2)k_2) |h_1|^{d/2+s -1/2} | h_2 |^{d/2-s+1/2} d h_1  d h_2 dk. \nonumber
		\end{align}
		Here the exchange of order of integration is justified by similar arguments as above. Combining the two outer integrals, we get that \eqref{formula::period-2nd} is equal to
		\begin{align*}
			\int_{(P^{(r)}_{2n})'(\BA) \backslash G_{2n}'(\BA)} &\int^{\ast}_{[G_{2n-r}]^{1,G_{2n-r}}}  W_r^{\psi} \left(\begin{pmatrix}
				m & \\ & 1_r
			\end{pmatrix}g,\phi \right) \mu_{s-1/2}(m) dm  \\
			& \qquad \cdot \Phi(e_2g) | \det g_1|^{s-n+(d+r-1)/2} | \det g_2|^{(d+r+1)/2-s-n} dg.
		\end{align*}
		In view of \eqref{formula::phi--varphi_P}, we get the disired formula by replacing $\phi$ by $\varphi_P$.
	\end{proof}

	We finish the proof of Theorem \ref{thm::1}. The arguments is similar to those in the proof of Proposition \ref{prop::explicit decomposition}. For each place $v$ of $F$, we use subscript to denote corresponding local objects. The crucial point here is that the regularized twisted period of square-integrable automorphic forms is factorizable. For cuspidal automorphic forms an explicit factorization is given by Proposition \ref{prop::explicit decomposition}. If an square-integrable automorphic form $\varphi$ is factorizable, then $\varphi_P$ is still factorizable  by Lemma \ref{lem::constant term--residue}, and so does $W_r^{\psi}(\varphi)$ by uniqueness of the Whittaker model. Then by Proposition \ref{prop::period formula} and induction, $\period(\varphi,s)$ is factorizable. Thus, we can choose local factors $\Theta^{\psi}$ such that 
	\begin{align*}
		\Theta^{\psi}(g,\varphi;s) = \prod_v \Theta^{\psi_v}(g_v,\varphi_v;s)
	\end{align*}
	for $\varphi = \otimes \varphi_v \in \speh(\sigma_v,d)$, where $\Theta^{\psi_v}(\mathbf{1},\varphi_{0,v};s) = 1$ for almost all spherical vectors $\varphi_{0,v} \in \speh(\sigma_v,d)$ used in the definition of the restricted tensor products. So, 
	\begin{align*}
		\hat{\Phi}(0) \period(\varphi,s) = r v_P \prod_v \int_{(P^{(r)}_{2n,v})' \backslash G_{2n,v}'} \Theta^{\psi_v}(g_v,\varphi_v;&s-1/2) \Phi_v(e_ng_{2,v}) \\
		& \cdot |g_{1,v}|^s|g_{2,v}|^{-s} |g_{2,v}| dg_{1,v} dg_{2,v}.
	\end{align*}
	By Proposition \ref{prop::Bump-F} and the proof of Proposition \ref{prop::period formula}, the local integrals and the infinite product converge absolutely.
	
	We define a function $\beta_v$ on $(P_{2n}')_v \backslash (G_{2n}')_v$ by
	\begin{align*}
		\beta_v(g, s ) = \int_{(P_n^{(r/2)})_v \times (P_{n-1}^{(r/2-1)})_v \backslash (G_n)_v \times (G_{n-1})_v} \Theta^{\psi_v}(\iota(g_1,g_2)g,\shskip\shskip &\varphi_v, s-1/2) \\
		& \cdot |g_1|^s |g_2|^{-s} dg_1dg_2.
	\end{align*}
	For $X \in F_v^{n-1}\oplus F_v^{\times}$, let $X_n \in F_v^{\times}$ be its $n$-th component and $n(X) \in G_{n,v}$ be the matrix with the last row $X$ and $n(X)_{ij} = \delta_{ij}$ for $i < n$. By Lemma \ref{lem::integration formula},
	\begin{align*}
		&   \int_{(P^{(r)}_{2n,v})' \backslash G_{2n,v}'} \Theta^{\psi_v}(g,\varphi_v;s-1/2) \Phi_v(e_ng_2) |g_1|^s|g_2|^{-s} |g_2| dg  \\
		=  &\int_{F_v^{n-1}\oplus F_v^{\times}} \Phi_v(X)\beta_v((\mathbf{1}_n;n(X)),s) |X_n|^{-s} dX
	\end{align*}
	If $\period(\varphi,s) = 0$, then there is a place $v$ such that 
	\begin{align*}
		\int_{F_v^{n-1}\oplus F_v^{\times}} \Phi_v(X) \beta_v ((\mathbf{1}_n;n(X)),s) |X_n|^{-s} dX = 0
	\end{align*}
	for all $\Phi_v \in \CS(F_v^n)$. Hence $\beta_v(\mathbf{1}_{2n},s) = 0$. Assume $\period(\varphi,s)$ is nonzero. Let $v_0$ be an arbitrary place of $F$. Fix $\Phi_v$ for $v \ne v_0$ and vary $\Phi_{v_0}$. We see that the local integral at $v_ 0$ is a constant of $\hat{\Phi}_{v_0}(0)$. Thus, $\beta_{v_0}((\mathbf{1}_n;n(X)),s)|X_n|^{-s}$ is a constant function which give $\beta_{v_0}(\mathbf{1}_{2n},s)$ when evaluated at $X = e_n$. Therefore the local integral equals to $\hat{\Phi}_{v_0}(0)\beta_{v_0}(\mathbf{1}_{2n},s)$, which completes the proof.

	\section{Distinguished residual spectrum}\label{sec::dist-spectrum}
	
	Recall that we have shown in Proposition \ref{prop::period formula} that, for $\varphi \in \speh(\sigma,d)$, $\period(\varphi,s)$ is analytic for all $s \in \BC$. We say that $\speh(\sigma,d)$ is $(G',\mu_{s_0})$-distinguished if there is a $\varphi \in \speh(\sigma,d)$ such that $\period(\varphi,s_0)$ is nonzero.

	By conjugating an Weyl element we see that $\varphi$ is $(G',\mu_{s})$-distingusihed if and only if $\varphi$ is $(G',\mu_{-s})$-distinguished.
	
	\begin{lem}
		Let $\sigma$ be an irreducible cuspidal automorphic representation of $G_r(\BA)$ with trivial central character. Assume that $\speh(\sigma,d)$ is $(G',\mu_s)$-distinguished. Then $r$ is an even integer and $\sigma$ is $(G_r',\mu_s)$-distinguished resp. $(G'_r,\mu_{s-1/2})$-distinguished when $d$ is even resp. $d$ is odd.
	\end{lem}
	\begin{proof}
		By \eqref{formula::period-2nd} or Theorem \ref{thm::1}, the regularized $(G',\mu_{s - 1/2})$-period on $\speh(\sigma,d-1)$ is realized as an inner integral of $(G',\mu_s)$-period on $\speh(\sigma,d)$. The lemma then follows by induction.
	\end{proof}
	
	We now show the converse. Let's retain notations from previous sections. Assume that $\speh(\sigma,d-1)$ is $(G'_{2n-r},\mu_{s-1/2})$-distinguished. Recall that we have shown above the regularized twisted period integral is factorizable. We have to show each local integral is nonzero.
	
	We switch to the local context and drop the subsript $v$ for simplicity. Let $F$ be a local field of characteristic zero. Let $\sigma$ be an irreducible unitary generic representation of $G_r(F)$. Note that $\speh(\sigma,d)$ is isomorphic to the unique irreducible subrepresentation of 
	\begin{align*}
		\Ind_P^G (\speh(\sigma,d-1)\boxtimes \sigma) \otimes \delta_{P}^{-1/(2r)}.
	\end{align*}
	By composing a Whittaker functional for the cuspidal part, we get a nonzero map
	\begin{align*}
		W_r^{\psi}: \speh(\sigma,d) \ra \ind_{P^{(r)}}^G \speh(\sigma,d-1)[(r-1)/2] \otimes \psi|_{N_r(F)}
	\end{align*}
	For $\varphi \in \speh(\sigma,d)$, the evaluation $W_r^{\psi}(e,\varphi)$ at $e$ gives a factorization of the function $W_r^{\psi}$ in the global situation. Suppose that $l$ is a nonzero $(G_{2n-r}',\mu_{1/2-s})$-invariant linear form on $\speh(\sigma,d-1)$. Set
	\begin{align*}
		\Theta^{\psi} (g,\varphi) = l(W_r^{\psi}(g,\varphi)).
	\end{align*}
	for $\varphi \in \speh(\sigma,d)$ and $g \in G_{2n}(F)$. Set 
	\begin{align}\label{formula::defn--beta-local form}
		\beta(\varphi)  =  \int_{P_n^{(r/2)} \times P_{n-1}^{(r/2-1) }  \backslash G_n \times G_{n-1}} \Theta^{\psi} \left( \iota(g_1,g_2),\varphi \right) |g_1|^s |g_2|^{-s} d(g_1,g_2).
	\end{align}
	It is easy to check that the integral is well-defined. If $l$ is the local invariant linear form coming from a factorization of the global period, then the integral in the right hand side of \eqref{formula::defn--beta-local form} is absolutely convergent as shown in the proof of Proposition \ref{prop::period formula}. Also in this case $\beta$ is nothing but the local integral, up to a nonzero constant, in a factorization of the global period. Thus, it remains to prove the following lemma.
	\begin{lem}
		Let the notations be as above. The linear form $\beta$ in \eqref{formula::defn--beta-local form} is nonzero on $\speh(\sigma,d)$.
	\end{lem}
	\begin{proof}
		We assume on the contrary that $\beta(\varphi) = 0$ for all $\varphi \in \speh(\sigma,d)$. For $\varphi \in \speh(\sigma,d)$, $ 0 \leqslant l \leqslant r/2$ and $1 \leqslant k \leqslant r/2$, define 
		\begin{align*}
			B_{l,k} ( g, \varphi)  =  \int_{P_{n-l}^{(r/2-l)} \times P_{n-k}^{(r/2-k) }  \backslash G_{n-l} \times G_{n-k}} \Theta^{\psi} \left( \left(   \begin{pmatrix}
				h_1 &  \\   & I_l
			\end{pmatrix}  ; \begin{pmatrix}
				h_2  &  \\  & I_k
			\end{pmatrix}  \right) g,\varphi  \right)       \\
			\cdot  |h_1|^{s-l}  |h_2|^{-s-k+1}  d(h_1,h_2).
		\end{align*}
		and $\beta_{l,k}(\varphi)  =  B_{l,k}(\mathbf{1}_{2n},\varphi)$. We will show recursively that $\beta_{r/2,r/2}(\varphi) =0$ which leads to a contradiction. By assumption, $\beta_{0,1}(\varphi) = \beta(\varphi) = 0$ for all $\varphi$. We then need only to show that $\beta_{l,k}(\varphi) \equiv 0$ implies $\beta_{l+1,k}(\varphi) \equiv 0$ when $l+1 = k$, and that $\beta_{l,k}(\varphi) \equiv 0$ imples $\beta_{l,k+1}(\varphi) \equiv 0$ when $l = k$. Since the proof of these two cases are similar, we only prove the first implication. 
		
		Suppose that $\beta_{l,k}(\varphi) = 0$ for all $\varphi$ with $l + 1 =k$. For $u = (u_1,u_2,\cdots,u_{n-l}) \in F^{n-l}$,  set 
		\begin{align*}
			\tilde{u} = (u_1,0,u_2,0,\cdots,0,u_{n-l}) \in F^{2n-k-l}	
		\end{align*}
		For $\Phi \in \CS(F^{n-l})$, define $\pi(\Phi)\varphi \in \pi$ by
		\begin{align*}
			\pi(\Phi)\varphi = \int_{F^{n-l}} \Phi(u) \pi  \begin{pmatrix}
				I_{2n-k-l} & {}^t\tilde{u} &  \\   & 1 &  \\ & & I_{k+l-1}
			\end{pmatrix}  \varphi du.
		\end{align*}
		Since
		\begin{align*}
			\iota \left(  \begin{pmatrix}
				h_1 & \\  & I_l
			\end{pmatrix} , \begin{pmatrix}
				h_2 &  \\ & I_k
			\end{pmatrix} \right)   \begin{pmatrix}
				I_{2n-k-l} & {}^t\tilde{u} &  \\   & 1 &  \\ & & I_{k+l-1}
			\end{pmatrix}     \\
			= \begin{pmatrix}
				I_{2n-k-l} & {}^t\widetilde{u{}^th} &  \\   & 1 &  \\ & & I_{k+l-1}
			\end{pmatrix}  \iota \left(  \begin{pmatrix}
				h_1 & \\  & I_l
			\end{pmatrix} , \begin{pmatrix}
				h_2 &  \\ & I_k
			\end{pmatrix} \right),
		\end{align*}
		we have
		\begin{align*}
			&	\Theta^{\psi} \left( \iota \left(  \begin{pmatrix}
				h_1 & \\  & I_l
			\end{pmatrix} , \begin{pmatrix}
				h_2 &  \\ & I_k
			\end{pmatrix} \right) ,\pi(\Phi) \varphi\right)  \\ & = \hat{\Phi}(e_{n-l}h_1) \Theta^{\psi} \left( \iota \left(  \begin{pmatrix}
				h_1 & \\  & I_l
			\end{pmatrix} , \begin{pmatrix}
				h_2 &  \\ & I_k
			\end{pmatrix} \right) , \varphi\right)
		\end{align*}
		By our assumption, $\beta_{l,k}(\pi(\Phi)\varphi) = 0$ for all $\Phi \in \CS(F^{n-l})$. By definition, $\beta_{l,k}(\pi(\Phi))$ equals to 
		\begin{align*}
			\int_{P_{n-l}^{(r/2-l)} \times P_{n-k}^{(r/2-k) }  \backslash G_{n-l} \times G_{n-k}} \hat{\Phi} (e_{n-l}h_1) \Theta^{\psi} &\left( \left(   \begin{pmatrix}
				h_1 &  \\   & I_l
			\end{pmatrix}  , \begin{pmatrix}
				h_2  &  \\  & I_k
			\end{pmatrix}  \right) ,\varphi  \right)   
			\\ &|h_1|^{s-l}  |h_2|^{-s-k+1}  d(h_1,h_2),
		\end{align*}
		which, by Lemma \ref{lem::integration formula}, equals to
		\begin{align*}
			\int_{\Xi^{n-l}} \hat{\Phi}(X) |X_{n-l}|^{s-l} B_{l+1,k}(n'(X),\varphi) dX.
		\end{align*}
		Here $n'(X)$ is the matrix obtained by replacing the $2n-k-l$ row of $\mathbf{1}_{2n}$ by 
		$$ (X_1,0,X_2,0,\cdots,0,X_{n-l},0,\cdots,0).$$
		Since $\hat{\Phi}$ is arbitary in $\CS(F^{n-l})$, we get $\beta_{l+1,k}(\varphi) = 0$.
		
	   The modifications needed for the archimedean case are exactly the same with those in \cite[Lemma 7.2]{Yamana-Residue-Periods-JFA}, so we omit the details.
	\end{proof}

	\section{Vanishing of periods when $p \ne q$}\label{sec:: case-p not q}
	
	In this short section let $G = G_n$ and $G' = M_{(p,q)}$, the standard Levi subgroup of $G_n$ associated to $(p,q)$. We assume $p > q$. The goal is to show the vanishing of the regularized (twisted) period of an square-integrable automorphic form on $G(\BA)$.
	
	For $s \in \BC$, we define a character $\xi_s$ on $G'(\BA)$ by
	\begin{align*}
		\xi_s \begin{pmatrix}
			g_1 &  \\
			   & g_2
		\end{pmatrix}  = |g_1|^{qs} |g_2|^{-ps}
	\end{align*}
	\begin{thm}
		Let $n = dr$ and $\sigma$ be an irreducible cuspidal automorphic representation of $G_r(\BA)$ with trivial central character. Let  $\varphi \in \speh(\sigma,d)$. Then, for all but a finite number of $s$, $\varphi$ is $\xi_s$-regular, and we have
		\begin{align*}
			\int^{\ast}_{[G']^{1,G}}\varphi (g) \xi_s (g) dg =  0.
		\end{align*}
	\end{thm}
	
	It is well known that the twisted period of a cusp form over $G'$ vanishes \cite{Ash-Ginzburg-Rallis-vanishing-pair},\cite[Proposition 2.1]{Jacquet-Friedberg-Crell-LinearPeriods}. Their proofs made use of the Fourier expansion along unipotent subgroups and worked for any cuspidal functions. The generalization to square-integrable automorphic forms here is a simple consequence of the local vanishing result Proposition \ref{prop::vanishing--unequal case} once we can show that the regularized integral is well-defined, as what we did in Lemma \ref{lem::vanishing--phi*E}. The computations are similar to those in Lemma \ref{lem::phi*E--regular}. We only point out the following fact. Let 
	\begin{align*}
		\xi = (\underbrace{q,\cdots,q}_p,\underbrace{-p,\cdots,-p}_q) \in \mathfrak{a}_{G'}^G.
	\end{align*}
	For $P \in \CP^{G,\textup{max}}(P'_0)$, $\langle \xi,\mathfrak{a}_P^G \rangle  = 0 $ implies 
	\begin{align*}
		\langle \underline{\rho}_P, t_P^G \rangle   \leqslant 0
	\end{align*}
	with $t_P^G$ defined in \eqref{formula::defn-t_P^Q}. We remark that the implication does not hold if $P$ is not maximal.
	
	\section*{Acknowledgements}
	
	The author thanks Nadir Matringe for his encouragement. He also would like to thank Atsushi Ichino and Shunsuke Yamana for answering quesions about their paper in email correspondences. This work is supported by the National Natural Science Foundation of China (No.12001191).

	\def\cprime{$'$} \def\cprime{$'$}

	
\end{document}